\documentclass[12pt]{amsart}
\usepackage{amsmath}
\usepackage{amssymb}
\usepackage{latexsym}

\newcommand{\Zint}{\mathbb {Z}}    
     
\newcommand{\Cplx}{\mathbb {C}}     
\newcommand{\halmos}{\rule{5pt}{5pt}}

\numberwithin{equation}{section}

\newtheorem{prop}{\bf Proposition}[section]
\newtheorem{thm}[prop]{\bf Theorem}

\newtheorem{cor}[prop]{\bf Corollary}

\begin{document}

\title[Integral representation of solutions]
{Integral representation of solutions to Fuchsian system and Heun's equation}
\author{Kouichi Takemura}
\dedicatory{Dedicated to Professor Masaki Kashiwara on his sixtieth birthday}
\address{Department of Mathematical Sciences, Yokohama City University, 22-2 Seto, Kanazawa-ku, Yokohama 236-0027, Japan.}
\email{takemura@yokohama-cu.ac.jp}
\begin{abstract}
We obtain integral representations of solutions to special cases of the Fuchsian system of differential equations and Heun's differential equation.
In particular, we calculate the monodromy of solutions to the Fuchsian equation that corresponds to Picard's solution of the sixth Painlev\'e equation, and to Heun's equation.
\end{abstract}

\maketitle

\section{Introduction}

The Fuchsian differential equation is a linear differential equation whose singularities are all regular.
It frequently appears in a range of problems in mathematics and physics.  
For example, the famous Gauss hypergeometric differential equation is a canonical form of the second-order Fuchsian differential equation with three singularities on the Riemann sphere $\Cplx \cup \{ \infty \}$.
Global properties of solutions, i.e., the monodromy, often play decisive roles in the applications of these equations in physics and other areas of mathematics.

Heun's differential equation is a canonical form of a second-order Fuchsian equation with four singularities, which is given by
\begin{equation}
\frac{d^2y}{dz^2} + \left( \frac{\gamma}{z}+\frac{\delta }{z-1}+\frac{\epsilon}{z-t}\right) \frac{dy}{dz} +\frac{\alpha \beta z -q}{z(z-1)(z-t)} y=0,
\label{Heun}
\end{equation}
with the condition 
\begin{equation}
\gamma +\delta +\epsilon =\alpha +\beta +1.
\label{Heuncond}
\end{equation}
Several approaches for analyzing Heun's equation are known: including the Heun polynomial (\cite{Ron}), Heun function (\cite{Ron}),
perturbation from the hypergeometric equation (\cite{Tak2}) and finite-gap integration (\cite{TV,GW,Smi,Tak3}).
 Finite-gap integration is applicable for the case $\gamma ,\delta ,\epsilon , \alpha - \beta \in \Zint +1/2$, $t \in \Cplx \setminus \{0,1 \} $ and all $q$, and results on the integral representation of solutions (\cite{Tak1}), the Bethe Ansatz (\cite{Tak1}), the Hermite-Krichever Ansatz (\cite{BE,Tak4}), the monodromy formulae by hyperelliptic integrals (\cite{Tak3}), the hyperelliptic-to-elliptic reduction formulae (\cite{Tak4}) and relationships with the Darboux transformation (\cite{Tak5}) have been obtained.
In this paper, we obtain integral formulae of solutions for the case $\gamma ,\delta ,\epsilon , \alpha +1/2,  \beta +1/2 \in \Zint $, $t \in \Cplx \setminus \{0,1 \} $ and all $q$, which then facilitates a calculation of  the monodromy.

To obtain these formulae, we need to consider a Fuchsian system of differential equations with four singularities $0,1,t,\infty $,
\begin{equation}
\frac{dY}{dz}=\left( \frac{A_0}{z}+\frac{A_1}{z-1}+\frac{A_t}{z-t} \right) Y,
\label{eq:dYdzAzY00}
\end{equation}
where $A_0$, $A_1$, $A_t$ are $2 \times 2$ matrices with constant elements. 
We consider the case that $\det A_0=\det A_1=\det A_t=0$, and $A_0+A_1+A_t=-\mbox{diag} (\kappa _1, \kappa _2)$ is a diagonal matrix.
Let $\theta _i$ $(i=0,1,t)$ denote the eigenvalues of $A_i$ other than $0$, and $\theta _{\infty } = \kappa _1 -\kappa _2$.
Under some assumptions the sixth Painlev\'e system is obtained by the monodromy preserving deformation.
Here the sixth Painlev\'e system is defined by
\begin{equation}
\frac{d\lambda }{dt} =\frac{\partial H_{VI}}{\partial \mu}, \quad \quad
\frac{d\mu }{dt} =-\frac{\partial H_{VI}}{\partial \lambda} ,
\label{eq:Psys}
\end{equation}
with the Hamiltonian 
\begin{align}
H_{VI} = & \frac{1}{t(t-1)} \left\{ \lambda (\lambda -1) (\lambda -t) \mu^2 \right. \label{eq:P6} \\
& \left. -\left\{ \theta _0 (\lambda -1) (\lambda -t)+\theta _1 \lambda (\lambda -t) +(\theta _t -1) \lambda (\lambda -1) \right\} \mu +\kappa _1(\kappa _2 +1) (\lambda -t)\right\} .\nonumber
\end{align}
By eliminating $\mu $ in Eq.(\ref{eq:Psys}), we obtain the sixth Painlev\'e equation for $\lambda $,
\begin{align}
\frac{d^2\lambda }{dt^2} = & \frac{1}{2} \left( \frac{1}{\lambda }+\frac{1}{\lambda -1}+\frac{1}{\lambda -t} \right) \left( \frac{d\lambda }{dt} \right) ^2 -\left( \frac {1}{t} +\frac {1}{t-1} +\frac {1}{\lambda -t} \right)\frac{d\lambda }{dt} \label{eq:P6eqn} \\
& +\frac{\lambda (\lambda -1)(\lambda -t)}{t^2(t-1)^2}\left\{ \frac{(1-\theta _{\infty})^2}{2} -\frac{\theta _{0}^2}{2}\frac{t}{\lambda ^2} +\frac{\theta  _{1}^2}{2}\frac{(t-1)}{(\lambda -1)^2} +\frac{(1-\theta _{t}^2)}{2}\frac{t(t-1)}{(\lambda -t)^2} \right\}, \nonumber
\end{align}
which is a non-linear ordinary differential equation of order two whose solutions do not have movable singularities other than poles.  
It is known that the sixth Painlev\'e systems have symmetry, and the action of the symmetry is called the Okamoto-B\"acklund transformation.
The sixth Painlev\'e system has two-parameter solutions for the case $\theta_0= \theta_1= \theta_t= 1-\theta_{\infty }=0$, which are called Picard's solution.
By the Okamoto-B\"acklund transformation of the sixth Painlev\'e system, Picard's solutions are transformed to the solutions for the case $(\theta _0 , \theta _1 , \theta _t, 1-\theta _{\infty} ) \in O_1 \cup O_2$, where
\begin{align}
& O_1= \left\{ (\theta _0 , \theta _1 , \theta _t, 1- \theta _{\infty} ) | 
\theta _0 , \theta _1 , \theta _t, 1- \theta _{\infty} \in \Zint +\frac{1}{2} \right\}, \\
& O_2 = \left \{(\theta _0 , \theta _1 , \theta _t, 1- \theta _{\infty} ) \left| 
\begin{array}{ll}
\theta _0 , \theta _1 , \theta _t, 1- \theta _{\infty}  \in \Zint \\
\theta _0 + \theta _1 + \theta _t + 1- \theta _{\infty}  \in 2 \Zint 
\end{array}
\right. \right\}. \nonumber 
\end{align}
For the case $(\theta _0 , \theta _1 , \theta _t, 1-\theta _{\infty} ) \in O_1$, solutions of the Fuchsian system (Eq.(\ref{eq:dYdzAzY00})) are expressed in the form of the Hermite-Krichever Ansatz, which is a consequence of results presented in \cite{TakP}.
In the present study, we investigate solutions of the Fuchsian system for the case $(\theta _0 , \theta _1 , \theta _t, 1-\theta _{\infty} ) \in O_2$.
 
These solutions will be shown to have integral representations whose integrands are functions in the form of the Hermite-Krichever Ansats.
In particular we obtain explicit solutions for the case $(\theta _0 , \theta _1 , \theta _t, 1-\theta _{\infty} ) =(0,0,0,0)$, and we can calculate the monodromy explicitly.
By considering the monodromy preserving deformation directly, we recover  Picard's solution of the sixth Painlev\'e equation.

The integral representions of solutions to the Fuchsian system follow from the results by Dettweiler-Reiter \cite{DR1,DR2} and Filipuk \cite{Fil} on the middle convolution (see section \ref{sec:MC}).
By considering special cases, we obtain integral formulae of solutions to Heun's equation for the case $\gamma ,\delta ,\epsilon , \alpha +1/2, \beta +1/2 \in \Zint $, $t \in \Cplx \setminus \{0,1 \} $ and all $q$, which are then available for calculating the monodromy.
For the case $\gamma =\delta =\epsilon =1, \alpha =3/2, \beta =1/2 $, we have explicit representations of the integral, and so we obtain explicit representations of the monodromy. 

This paper is organized as follows:
In section \ref{sec:FS}, we introduce notation for the Fuchsian system with four singularities.
In section \ref{sec:MC}, we review results on the middle convolution due to Dettweiler-Reiter and Filipuk, and combine their results.
In section \ref{sec:HKA}, we recall the Hermite-Krichever Ansatz.
In section \ref{sec:intrepFs}, we obtain integral representations of solutions to the Fuchsian system with four singularities for the case $(\theta _0 , \theta _1 , \theta _t, 1-\theta _{\infty} ) \in O_2$, whose integrands are functions in the form of  the Hermite-Krichever Ansats.
In section \ref{sec:intrepFs0000}, we have explicit representations of solutions for the case $(\theta _0 , \theta _1 , \theta _t, 1-\theta _{\infty} ) =(0,0,0,0)$, and we calculate the monodromy explicitly.
Furthermore, by considering the monodromy preserving deformation directly, we recover Picard's solution of the sixth Painlev\'e equation.
In section \ref{sec:intrepHeun}, we obtain integral formulae of solutions to Heun's equation for the case $\gamma ,\delta ,\epsilon , \alpha - \beta -1/2 \in \Zint $, $t \in \Cplx \setminus \{0,1 \} $ and all $q$, including the case $\gamma =\delta =\epsilon =1, \alpha =3/2, \beta =1/2 $.

\section{Fuchsian system with four singularities} \label{sec:FS}

We consider a system of ordinary differential equations,
\begin{equation}
\frac{dY}{dz}=A(z)Y, \quad A(z)=\frac{A_0}{z}+\frac{A_1}{z-1}+\frac{A_t}{z-t} =
\left(
\begin{array}{ll}
a_{11}(z) & a_{12}(z) \\
a_{21}(z) & a_{22}(z) 
\end{array}
\right) ,
\label{eq:dYdzAzY}
\end{equation}
where $Y= {}^t (y_1(z), y_2(z)) $, $t\neq 0,1$ , $A_0$, $A_1$, $A_t$ are $2 \times 2$ matrices with constant elements.
Then Eq.(\ref{eq:dYdzAzY}) is Fuchsian, i.e., any singularities on the Riemann sphere $\Cplx \cup \{ \infty \} $ are regular, and it may have regular singularities at $z=0,1,t,\infty$ on this sphere.
Set
\begin{align}
& A_0= \left(
\begin{array}{ll}
u_0+\theta _0  & -w_0 \\
u_0(u_0+\theta _0)/ w_0  & -u_0
\end{array}
\right) , \quad 
A_1= \left(
\begin{array}{ll}
u_1+\theta _1  & -w_1 \\
u_1(u_1+\theta _1)/w_1 & -u_1
\end{array}
\right) , \label{eq:A0A1AtP} \\
& A_t= \left(
\begin{array}{ll}
u_t+\theta _t  & -w_t \\
u_t(u_t+\theta _t)/w_t  & -u_t
\end{array}
\right) , \quad  \nonumber
\end{align}
where $u_0, w_0, u_1, w_1, u_t, w_t$ are defined by
\begin{align}
& w_0 = \frac{k\lambda}{t}, \quad  w_1= -\frac{k(\lambda-1)}{t-1}, \quad w_t = \frac{k(\lambda-t)}{t(t-1)} , \label{eq:wugen}\\
& u_0=-\theta _0 +\frac{\lambda}{t\theta _{\infty}} [ \lambda(\lambda-1)(\lambda-t)\mu^2 +\{ 2\kappa _1 (\lambda-1)(\lambda-t)-\theta _1(\lambda-t) \nonumber \\
& \quad \quad \quad \quad \quad \quad \quad \quad -t\theta _t(\lambda-1) \} \mu +\kappa _1 \{\kappa _1(\lambda-t-1)-\theta _1-t\theta _t\} ] ,\nonumber \\
& u_1 =-\theta _1 -\frac{\lambda-1}{(t-1)\theta _{\infty}}  [ \lambda(\lambda-1)(\lambda-t)\mu^2 + \{ 2\kappa _1 (\lambda-1)(\lambda-t)+(\theta _{\infty}-\theta _1 )(\lambda-t) \nonumber\\
& \quad \quad \quad \quad \quad \quad \quad \quad -t\theta _t (\lambda-1)\}\mu +\kappa _1 \{ \kappa _1 (\lambda-t+1)+\theta _0-(t-1)\theta _t\} ] ,\nonumber \\
& u_t= -\theta _t+\frac{\lambda-t}{t(t-1)\theta _{\infty}}  [ \lambda(\lambda-1)(\lambda-t)\mu^2 + \{ 2\kappa _1 (\lambda-1)(\lambda-t) -\theta _1 (\lambda-t) \nonumber\\
& \quad \quad \quad \quad \quad \quad \quad \quad +t(\theta _{\infty}-\theta _t )(\lambda-1)\} \mu +\kappa _1 \{ \kappa _1 (\lambda-t+1)+\theta _0+(t-1)(\theta _{\infty}-\theta _t)\} ] ,\nonumber 
\end{align}
and $\kappa _1= (\theta _{\infty } -\theta _0 -\theta _1 -\theta _t)/2$, $\kappa _2= -(\theta _{\infty } +\theta _0 +\theta _1 +\theta _t)/2$.
Note that the eigenvalues of $A_i$ $(i=0,1,t)$ are $0$ and $\theta _i$.
Set $A_{\infty }= -(A_0 +A_1 +A_t)$. Then 
\begin{align}
A_{\infty }= \left(
\begin{array}{cc}
\kappa _1 & 0 \\
0 & \kappa _2 
\end{array}
\right) . \label{def:Ainf}
\end{align}
We denote the Fuchsian sysytem (Eq.(\ref{eq:dYdzAzY})) with Eqs.(\ref{eq:A0A1AtP}, \ref{eq:wugen}) by $D_Y(\theta _0, \theta _1, \theta _t, \theta _{\infty}; \lambda ,\mu ;k )$.
By eliminating $y_2(z)$ in Eq.(\ref{eq:dYdzAzY}), we have a second-order linear differential equation,
\begin{align}
& \frac{d^2y_1(z)}{dz^2} + \left( \frac{1-\theta _0}{z}+\frac{1-\theta _1}{z-1}+\frac{1-\theta _t}{z-t}-\frac{1}{z-\lambda}  \right)  \frac{dy_1(z)}{dz} \label{eq:linP6} \\
&  \quad \quad + \left( \frac{\kappa _1(\kappa _2 +1)}{z(z-1)}+\frac{\lambda (\lambda -1)\mu}{z(z-1)(z-\lambda)}-\frac{t (t -1)H}{z(z-1)(z-t)}  \right) y_1(z)=0, \nonumber \\
& H=\frac{1}{t(t-1)}[ \lambda (\lambda -1) (\lambda -t)\mu ^2 -\{ \theta _0  (\lambda -1) (\lambda -t)+\theta _1  \lambda (\lambda -t) \nonumber \\
& \quad \quad \quad \quad \quad \quad \quad \quad +(\theta _t -1) \lambda (\lambda -1)\} \mu +\kappa _1 (\kappa _2 +1) (\lambda -t)], \nonumber
\end{align}
which we denote by $D_{y_1}(\theta _0, \theta _1, \theta _t, \theta _{\infty}; \lambda ,\mu )$.
This equation has regular singularities at $z=0,1,t,\lambda ,\infty$. The exponents of the singularity $z=\lambda $ are $0,2$, and this singularity is apparent (i.e. non-logarithmic).
Note that the sixth Painlev\'e system
\begin{align}
\frac{d\lambda }{dt} =\frac{\partial H}{\partial \mu}, \quad \frac{d\mu }{dt} =-\frac{\partial H}{\partial \lambda}
\label{eq:P6sys}
\end{align}
describes the condition for the monodromy preserving deformation of Eq.(\ref{eq:dYdzAzY}) with respect to the variable $t$.

It is known that the sixth Painlev\'e system has symmetry of the extended affine Weyl group of type $F_{4}^{(1)}$ (\cite{Oka}), which is called the Okamoto-B\"acklund transformation.
In particular the sixth Painlev\'e system is invariant under Okamoto's transformation $s_2$ defined by
\begin{align}
 s_2: \; & \theta _0 \rightarrow \kappa _1+ \theta _0, \; \; \theta _1 \rightarrow \kappa _1+ \theta _1, \; \; \theta _t \rightarrow \kappa _1+ \theta _t, \; \; \theta _{\infty } \rightarrow -\kappa _2, \\ 
& \lambda \rightarrow \lambda  +\kappa _1/\mu ,\; \; \mu \rightarrow \mu, \; \; t \rightarrow t. \nonumber
\end{align}
Note that $s_2$ is involutive, i.e., $(s_2)^2=1$.

\section{Middle convolution} \label{sec:MC}

Dettweiler and Reiter \cite{DR1,DR2} gave an algebraic analogue of Katz' middle convolution functor, and Filipuk \cite{Fil} applied them for the Fuchsian system with four singularities.
We review and combine these authors' results for the present setting.  
Note that the results of Dettweiler and Reiter are valid for Fuchsian equations of an arbitrary size and an arbitrary number of singular points.
Let $A_0$, $A_1$, $A_t$ be matrices in $\Cplx ^{2\times 2}$.
For $\nu \in \Cplx$, we define the convolution matrices $B_0, B_1, B_t \in \Cplx ^{6\times 6}$ as follows:
\begin{align}
& B_0=
\left(
\begin{array}{ccc}
A_0 +\nu & A_1 & A_t \\
0 & 0 & 0 \\
0 & 0 & 0 
\end{array}
\right) , \quad
B_1=
\left(
\begin{array}{ccc}
0 & 0 & 0 \\
A_0 & A_1 +\nu & A_t \\
0 & 0 & 0 
\end{array}
\right) , \quad  \label{eq:Bdef} \\
& B_t=
\left(
\begin{array}{ccc}
0 & 0 & 0 \\
0 & 0 & 0 \\
A_0 & A_1 & A_t +\nu
\end{array}
\right) . \nonumber
\end{align}
We consider the following differential equation:
\begin{equation}
\frac{dU}{dz}=\left( \frac{B_0}{z}+\frac{B_1}{z-1}+\frac{B_t}{z-t} \right) U, \quad U \in \Cplx ^6
\label{eq:dYdzBzU}
\end{equation}
We fix a base point $o \in \Cplx \setminus \{0,1,t\}$.
Let $\alpha _i$ $(i=0,1,t,\infty )$ be a cycle turning the point $w=i$ anti-clockwise whose base point is $o$.
Let $z \in \Cplx \setminus \{0,1,t\}$
and $\alpha _{z}$ be a cycle turning the point $w=z$ anti-clockwise.
Let $[\alpha , \beta] = \alpha ^{-1} \beta ^{-1}\alpha \beta $ be the Pochhammer contour.
\begin{prop} $($\cite{DR2}$)$ \label{prop:DRintegrepr}
Assume that $Y= {}^t (y_1(z), y_2(z))$ is a solution to the differential equation
\begin{equation}
\frac{dY}{dz}=\left( \frac{A_0}{z}+\frac{A_1}{z-1}+\frac{A_t}{z-t} \right) Y.
\label{eq:dYdzAzY0}
\end{equation}
For $i \in \{ 0,1,t,\infty \}$, the function
\begin{equation}
U = \left(
\begin{array}{l}
\int _{[\alpha _{z} ,\alpha _i]}w^{-1}y_1(w) (z-w)^{\nu } dw \\
\int _{[\alpha _{z} ,\alpha _i]}w^{-1}y_2(w) (z-w)^{\nu } dw \\
\int _{[\alpha _{z} ,\alpha _i]}(w-1)^{-1}y_1(w) (z-w)^{\nu } dw \\
\int _{[\alpha _{z} ,\alpha _i]}(w-1)^{-1}y_2(w) (z-w)^{\nu } dw \\
\int _{[\alpha _{z} ,\alpha _i]}(w-t)^{-1}y_1(w) (z-w)^{\nu } dw \\
\int _{[\alpha _{z} ,\alpha _i]}(w-t)^{-1}y_2(w) (z-w)^{\nu } dw
\end{array}
\right) , \label{eq:integrepU}
\end{equation}
satisfies differential equation (\ref{eq:dYdzBzU}).
\end{prop}
\begin{proof}
It follows from a straightforward calculation that the function
\begin{equation}
U = \left(
\begin{array}{l}
z^{-1}y_1(z)\\
z^{-1}y_2(z)\\
(z-1)^{-1} y_1(z)\\
(z-1)^{-1} y_2(z)\\
(z-t)^{-1} y_1(z)\\
(z-t)^{-1} y_2(z)
\end{array}
\right) ,
\end{equation}
is a solution of Eq.(\ref{eq:dYdzBzU}) for the case $\nu =-1$ (see \cite[Lemma 6.4]{DR2}).

It is shown in \cite[Lemma 6.2]{DR2} that if $U={}^t (u_1 (z) , u_2 (z) , \dots , u_6 (z) )$ is a solution of Eq.(\ref{eq:dYdzBzU}) for the case $\nu =\nu _1$, then the function 
\begin{equation}
\bar{U} = \left( 
\begin{array}{c}
\int _{[\alpha _{z} ,\alpha _i]} u_1 (w)(w-z)^{\nu_2 -1} dw\\
\int _{[\alpha _{z} ,\alpha _i]} u_2 (w)(w-z)^{\nu_2 -1} dw\\
\vdots \\
\int _{[\alpha _{z} ,\alpha _i]} u_6 (w)(w-z)^{\nu_2 -1} dw\\
\end{array}
\right) ,
\end{equation}
is a solution of Eq.(\ref{eq:dYdzBzU}) for the case $\nu =\nu _1 +\nu _2$.
By applying this result for the case $\nu_1=-1, \nu_2=\nu +1$, we obtain the proposition.
\end{proof}
We set 
\begin{align}
& {\mathcal L}_0= \left(
\begin{array}{c}
\mbox{Ker}(A_0) \\
0 \\
0 
\end{array}
\right) , \quad 
{\mathcal L}_1= \left(
\begin{array}{c}
0\\
\mbox{Ker}(A_1) \\
0 
\end{array}
\right) , \quad 
{\mathcal L}_t= \left(
\begin{array}{c}
0 \\
0 \\
\mbox{Ker}(A_t) 
\end{array}
\right) , \\
& {\mathcal L} ={\mathcal L}_0 \oplus {\mathcal L}_1 \oplus {\mathcal L}_t, \quad 
{\mathcal K} = \mbox{Ker}(B_0) \cap \mbox{Ker}(B_1) \cap  \mbox{Ker}(B_t) . \nonumber
\end{align}
We fix an isomorphism between $\Cplx ^6 /({\mathcal K}+ {\mathcal L})$ and $\Cplx ^m$ for some $m$.
A tuple of matrices $mc_{\nu } (A) =(\tilde{B}_0, \tilde{B}_1, \tilde{B}_t)$, where $\tilde{B}_k$ $(k=0,1,t)$ is induced by the action of $B_k$ on $\Cplx ^m \simeq \Cplx ^6 /({\mathcal K}+ {\mathcal L})$, is called an additive version of the middle convolution of $(A_0,A_1,A_t)$ with the parameter $\nu$.
Filipuk \cite{Fil} established that, if $\nu =\kappa _1$, then $\Cplx ^6 /({\mathcal K}+ {\mathcal L})$ is isomorphic to $\Cplx ^2$ and the isomonodromic deformation of the middle convolution system
\begin{equation}
\frac{d\tilde{Y}}{dz}=\left( \frac{\tilde{B}_0}{z}+\frac{\tilde{B}_1}{z-1}+\frac{\tilde{B}_t}{z-t} \right) \tilde{Y},
\end{equation}
gives the sixth Painleve equation for the parameters transformed by Okamoto's transformation $s_2$.
Note that Boalch \cite{Boa} obtained a geometric result on Okamoto's transformation earlier by finding an isomorphism between a $2 \times 2$ Fuchsian equation and a $3 \times 3$ Fuchsian equation, which would be related with Filipuk's result.

We now calculate explicitly the Fuchsian differential equation determined by the middle convolution that is required for our purpose, and which reproduces the result by Filipuk \cite{Fil}. 
Let $A_0$, $A_1$, $A_t$ be the matrices defined by Eq.(\ref{eq:A0A1AtP}). If $\nu = \kappa _1$, then the spaces ${\mathcal L}_0$, ${\mathcal L}_1$, ${\mathcal L}_t$, ${\mathcal K}$ are written as
\begin{align}
& {\mathcal L}_0= \Cplx \left(
\begin{array}{c}
w_0 \\
u_0+\theta _0 \\
0\\ 
0\\
0\\
0
\end{array}
\right) , \;
{\mathcal L}_1= \Cplx \left(
\begin{array}{c}
0\\
0\\
w_1 \\
u_1+\theta _1 \\
0\\ 
0
\end{array}
\right) , \;
{\mathcal L}_t= \Cplx \left(
\begin{array}{c}
0\\
0\\
0\\
0\\
w_t \\
u_t+\theta _t
\end{array}
\right) , \;
{\mathcal K}= \Cplx \left(
\begin{array}{c}
1\\
0\\
1\\
0\\
1 \\
0
\end{array}
\right) .
\end{align}
Set
\begin{align}
& S= \left(
\begin{array}{cccccc}
0 & 0 & 1 & w_0 & 0 & 0\\
0 & 0 & 0 &u_0+\theta _0& 0 & 0 \\
s_{31} & s_{32} & 1 &0 & w_1 & 0 \\ 
0 & 0 & 0 &0 & u_1+\theta _1& 0 \\
s_{51} & s_{52} & 1 &0& 0 &  w_t \\
0 & 0 & 0 &0 & 0 & u_t+\theta _t
\end{array}
\right) , \;
\begin{array}{l}
s_{31}= \frac{t(t-1)\mu w_0 w_1 u_t}{k^2 \kappa _2 w_t(u_0+\theta _0)(u_1+\theta _1)}, \\
s_{51}= \frac{t(1-t)\mu w_0 w_t u_1}{k^2 \kappa _2 w_t(u_0+\theta _0)(u_t+\theta _t)}, \\
s_{32}= \frac{1}{\theta _{\infty }} \left( \frac{w_0}{u_0+\theta _0}- \frac{w_1}{u_1+\theta _1} \right) ,\\
s_{52}= \frac{1}{\theta _{\infty }} \left( \frac{w_0}{u_0+\theta _0}- \frac{w_t}{u_t+\theta _t} \right) ,
\end{array}
\end{align}
and $\tilde{U}= S^{-1} U$,  where $U$ is a solution to Eq.(\ref{eq:dYdzBzU}).
Then $\det U =k \lambda (\lambda -1)(\lambda -t) \mu/(t(1-t)\theta _{\infty})$ and $\tilde{U}$ satisfies
\begin{equation}
\frac{d\tilde{U}}{dz} = \left(
\begin{array}{cccccc}
b_{11}(z) & b_{12}(z) & 0 & 0 & 0 & 0\\
b_{21}(z) & b_{22}(z) & 0 & 0 & 0 & 0\\
b_{31}(z) & b_{32}(z) & 0 & 0 & 0 & 0\\
0 & b_{42}(z) & 0 & \frac{\kappa _1}{z} & 0 & 0\\
0 & b_{52}(z) & 0 & 0 & \frac{\kappa _1}{z-1} & 0\\
0 & b_{62}(z) & 0 & 0 & 0 & \frac{\kappa _1}{z-t}
\end{array}
\right) \tilde{U},
\end{equation}
where $b_{i1}(z)$ $(i=1,2,3)$ and $b_{i2}(z)$ $(i=1,\dots ,6)$ are rational functions.
Write $\tilde{U}= {}^t (\tilde{u}_1(z), \tilde{u}_2(z), \dots , \tilde{u}_6(z))$ and set $\tilde{y}_1(z) = \tilde{u}_1(z)$, $\tilde{y}_2(z) = \tilde{u}_2(z)$ and $\tilde{Y}= {}^t (\tilde{y}_1(z), \tilde{y}_2(z))$.
Then we have
\begin{align}
& \frac{d\tilde{Y}}{dz}=
\left(
\begin{array}{ll}
b_{11}(z) & b_{12}(z) \\
b_{21}(z) & b_{22}(z) 
\end{array}
\right) \tilde{Y}.
\label{eq:dtYdzBztY0} 
\end{align}
The elements $b_{11}(z)$, $b_{12}(z)$, $b_{21}(z)$, $b_{22}(z)$ are calculated explicitly and Eq.(\ref{eq:dtYdzBztY0})  coincides with the Fuchsian differential equation $D_Y(\tilde{\theta }_0, \tilde{\theta }_1, \tilde{\theta }_t, \tilde{\theta }_{\infty}; \tilde{\lambda },\tilde{\mu };\tilde{k} )$ (see Eq.(\ref{eq:dYdzAzY})), where 
\begin{align}
& \tilde{\theta }_0 = \frac{\theta _0 -\theta _1 -\theta _t+\theta _{\infty}}{2}, \quad \tilde{\theta }_1 = \frac{-\theta _0 +\theta _1 -\theta _t+\theta _{\infty}}{2}, \quad \tilde{\theta }_t = \frac{-\theta _0 -\theta _1 +\theta _t+\theta _{\infty}}{2}, \\
&  \tilde{\theta }_{\infty} = \frac{\theta _0+\theta _1 +\theta _t+\theta _{\infty}}{2}, \quad \tilde{\lambda} =\lambda + \kappa _1/\mu, \quad  \tilde{\mu} =\mu ,\quad \tilde{k}= k . \nonumber
\end{align}
The functions $\tilde{y}_1(z)$ and $\tilde{y}_2(z)$ are expressed as 
\begin{align}
& \tilde{y}_1(z) = \frac{\tilde{\lambda } (u_0+\theta _0)}{\lambda } u_1(z) - \frac{k\tilde{\lambda }}{t} u_2 (z) + \frac{(\tilde{\lambda } -1) (u_1+\theta _1)}{\lambda -1} u_3(z) \label{eq:y1tu16} \\
& \quad \quad \quad + \frac{k(\tilde{\lambda }-1)}{t-1 } u_4 (z)+ \frac{(\tilde{\lambda } -t)(u_t+\theta _t)}{\lambda -t} u_5(z) + \frac{k(\tilde{\lambda }-t)  }{t(1-t)} u_6 (z), \nonumber \\
& \tilde{y}_2(z) = \frac{\theta _{\infty}}{\kappa _2} \left( -\frac{t u_0 (u_0+\theta _0)}{k \lambda } u_1(z) + u_0 u_2 (z) + \frac{(t-1) u_1 (u_1 + \theta _1 )}{k (\lambda -1)} u_3(z) \right. \nonumber \\
& \left. \quad \quad \quad \quad \quad \quad + u_1 u_4 (z)+ \frac{t(1-t) u_t (u_t+\theta _t)}{k(\lambda -t)} u_5(z) +u_t u_6 (z) \right) . \nonumber
\end{align}
It follows from Proposition \ref{prop:DRintegrepr} that the function $U={}^t (u_1 (z) , u_2 (z) , \dots , u_6 (z) )$ given by Eq.(\ref{eq:integrepU}) is a solution to Eq.(\ref{eq:dYdzBzU}).  
Combining with the relations $y_2 (w)= (dy_1(w)/dw -a_{11}(w)y_1(w))/a_{12}(w)$, $y_1 (w)= (dy_2(w)/dw -a_{22}(w)y_2(w))/a_{21}(w)$ and Eq.(\ref{eq:y1tu16}),
the functions $\tilde{y}_1(z)$ and $\tilde{y}_2(z)$ are expressed as the integral in the following proposition by means of a straightforward calculation:
\begin{prop} \label{thm:zinterep}
Set $\kappa _1= (\theta _{\infty } -\theta _0 -\theta _1 -\theta _t)/2$ and $\kappa _2= -(\theta _{\infty } +\theta _0 +\theta _1 +\theta _t)/2$.
If $Y={}^t ( y_1(z), y_2(z))$ is a solution to the Fuchsian differential equation $D_Y(\theta _0, \theta _1, \theta _t,\theta _{\infty}; \lambda ,\mu ;k )$ (see Eq.(\ref{eq:dYdzAzY})), then the function $\tilde{Y} = {}^t (\tilde{y}_1 (z) , \tilde{y}_2 (z) )$ defined by 
\begin{align}
& \tilde{y}_1 (z)= \int _{[\alpha _{z} ,\alpha _i]} \left\{ \kappa _1 y_1(w) + (w- \tilde{\lambda})\frac{dy_1(w)}{dw}  \right\} \frac{(z-w)^{\kappa _1 }}{w-\lambda } dw,  \label{eq:yt1zintrep} \\
& \tilde{y}_2 (z) = \frac{-\theta _{\infty}}{\kappa _2 }\int _{[\alpha _{z} ,\alpha _i]}  \frac{dy_2(w)}{dw}  (z-w)^{\kappa _1 } dw,  \nonumber
\end{align}
satisfies the Fuchsian differential equation $D_Y(\kappa _1 + \theta _0, \kappa _1 +\theta _1, \kappa _1 +\theta _t, -\kappa _2 ; \lambda +\kappa _1 /\mu ,\mu ;k )$ for $i \in \{ 0,1,t,\infty \}$.
\end{prop}
Therefore, if we know a solution to the differential equation $D_Y(\theta _0, \theta _1, \theta _t,\theta _{\infty}; \lambda ,\mu ;k )$, then we have integral representations of solutions to the Fuchsian differential equation $D_Y(\tilde{\theta }_0, \tilde{\theta }_1, \tilde{\theta }_t, \tilde{\theta }_{\infty}; \tilde{\lambda },\tilde{\mu };\tilde{k} )$ obtained by Okamoto's transformation $s_2$.
It can be shown that, if $\kappa _2 \neq 0$, $\kappa _1\not \in \Zint$ and $\theta _i \not \in \Zint $ for some $i \in \{ 0,1,t,\infty \}$, then the function $\tilde{y}_1 (z)$ is non-zero for generic $\lambda $ and $\mu$ (see \cite[Lemma 6.6]{DR2}). 
On the other hand, a solution to Eq.(\ref{eq:dYdzAzY}) for the case $\theta _0, \theta _1, \theta _t, \theta _{\infty } \in \Zint +\frac{1}{2}$ can be expressed in the form of  the Hermite-Krichever Ansats.
In the next section, we recall the Hermite-Krichever Ansats.

\section{Hermite-Krichever Ansatz} \label{sec:HKA}
We rewrite Eq.(\ref{eq:linP6}) in  elliptical form. Recall that Eq.(\ref{eq:linP6}) is written as
\begin{align}
&  \frac{d^2y_1(z)}{dz^2} + \left( \frac{1-\theta _0}{z}+\frac{1-\theta _1}{z-1}+\frac{1-\theta _t}{z-t}-\frac{1}{z-\lambda}  \right)  \frac{dy_1(z)}{dz} \label{eq:linP60} \\
&  \quad \quad + \left( \frac{\kappa _1(\kappa _2 +1)}{z(z-1)}+\frac{\lambda (\lambda -1)\mu}{z(z-1)(z-\lambda)}-\frac{t (t -1)H}{z(z-1)(z-t)}  \right) y_1(z)=0, \nonumber 
\end{align}
and $H$ is determined by
\begin{align}
& H=\frac{1}{t(t-1)}[ \lambda (\lambda -1) (\lambda -t)\mu ^2 -\{ \theta _0  (\lambda -1) (\lambda -t)+\theta _1  \lambda (\lambda -t) \label{eq:linP6H} \\
& \quad \quad \quad \quad \quad \quad \quad \quad +(\theta _t -1) \lambda (\lambda -1)\} \mu +\kappa _1 (\kappa _2 +1) (\lambda -t)]. \nonumber
\end{align}
Let $\wp (x)$ be the Weierstrass $\wp$-function with periods $(2\omega _1,2\omega _3)$, $\omega _0(=0)$, $\omega_1$, $\omega_2(=-\omega _1 -\omega _3)$, $\omega_3$ be half-periods and $e_i=\wp (\omega _i)$ $(i=1,2,3)$.
Set
\begin{equation}
z=\frac{\wp (x)-e_1}{e_2-e_1} ,\quad t=\frac{e_3-e_2}{e_1-e_2}, \quad \lambda =\frac{\wp (\delta )-e_1}{e_2-e_1}
\end{equation}
For $t \in \Cplx \setminus \{0,1 \}$, there exists a pair of periods  $(2\omega _1,2\omega _3)$ such that $ t=(\wp (\omega _3)-\wp (\omega _2))/(\wp (\omega _1)-\wp (\omega _2))$. The value $\delta $ is determined up to the sign $\pm $ and the periods $2\omega _1 \Zint \oplus 2\omega _3 \Zint $.
Set
\begin{align}
& \theta _0 =l_1 +1/2, \quad \theta _1 =l_2 +1/2, \quad \theta _t =l_3 +1/2, \quad \theta _{\infty} =-l_0 +1/2, \label{eq:kili} \\
& f(x)= y_1 (z) z^{-l_1/2} (z-1)^{-l_2/2} (z-t)^{-l_3/2}. \nonumber
\end{align}
Then Eq.(\ref{eq:linP60}) is transformed to 
\begin{align}
& \left( -\frac{d^2}{dx^2} + \frac{\wp ' (x)}{\wp (x) -\wp (\delta )} \frac{d}{dx} + \frac{\tilde{s}}{\wp (x) -\wp (\delta )} +\sum_{i=0}^3 l_i(l_i+1) \wp (x+\omega_i) +C \right) f (x)=0, \label{eq:Hg} \\
& \tilde{s}= -4(e_2-e_1)^2 \lambda (\lambda -1) (\lambda -t) \left( \mu -\frac{l_1}{2\lambda } -\frac{l_2}{2 (\lambda -1)}-\frac{l_3}{2 (\lambda -t)} \right) , \nonumber \\
& C=4(e_2-e_1)\left\{ \lambda (1-\lambda )\mu- t(1-t)H\right\} +(l_1+l_2+l_3+l_0+1)(l_1+l_2+l_3-l_0)e_3 \nonumber \\
&  \quad \quad  -2(l_1l_2e_3+l_2l_3e_1+l_3l_1e_2) + 2(l_1+l_2+l_3) ((e_2-e_1)\lambda +e_1) +\sum _{i=1}^3 l_i (l_i +2) e_i, \nonumber 
\end{align}
and Eq.(\ref{eq:linP6H}) is equivalent to the equality
\begin{align}
& C =4(e_2-e_1) \lambda (\lambda -1) (\lambda -t) \mu \left\{ \mu - \frac{l_1+\frac{1}{2}}{\lambda } - \frac{ l_2+\frac{1}{2}}{\lambda -1} - \frac{l_3+\frac{1}{2}}{\lambda -t} \right\}  +\sum _{i=1}^3 l_i (l_i +2) e_i \label{pgkrap} \\
& \quad +((e_2-e_1) \lambda +e_1) \{ (l_1+l_2+l_3+l_0+2)(l_1+l_2+l_3-l_0+1) -2\} \nonumber \\
& \quad -2(l_1l_2e_3+l_2l_3e_1+l_3l_1e_2) , \nonumber
\end{align}
which shows that the regular singularities $x=\pm \delta $ are apparent.
The sixth Painlev\'e equation (Eq.(\ref{eq:P6eqn})) for $\lambda (=(\wp (\delta )-e_1)/(e_2-e_1))$ also has an elliptical representation
\begin{equation}
\frac{d^2 \delta }{d \tau ^2} = -\frac{1}{4\pi ^2} \left\{ \frac{(1-\theta _{\infty})^2}{2} \wp ' \left(\delta  \right) + \frac{\theta _{0}^2}{2} \wp ' \left(\delta +\frac{1}{2} \right) + \frac{\theta _{1}^2}{2} \wp ' \left(\delta +\frac{\tau +1}{2} \right) +  \frac{\theta _{t}^2}{2} \wp ' \left(\delta +\frac{\tau }{2}\right) \right\}, \label{eq:P6ellip}
\end{equation}
where  $\omega _1=1/2$, $\omega _3=\tau /2$ and $\wp ' (z ) = (\partial /\partial z ) \wp (z)$ (see \cite{Man,Tks,TakHP}),
and it is related to the monodromy preserving deformation of Eq.(\ref{eq:Hg}) by the variable $\tau =\omega _3/\omega _1$.

We recall that a solution to Eq.(\ref{eq:Hg}) can be expressed in the form of the Hermite-Krichever Ansatz if $l_0 ,l_1,l_2, l_3 \in \Zint$.
Note that the condition $\theta _{\infty} , \theta _0, \theta _1 ,  \theta _t \in \Zint +1/2$ corresponds to the condition $l_0 ,l_1,l_2, l_3 \in \Zint$.
Set
\begin{equation}
\Phi _i(x,\alpha )= \frac{\sigma (x+\omega _i -\alpha ) }{ \sigma (x+\omega _i )} \exp (\zeta( \alpha )x), \quad \quad (i=0,1,2,3).
\label{Phii}
\end{equation}
\begin{prop} \label{thm:alpha} $($\cite{TakP}$)$
Set $\tilde{l} _0 = |l_0 +1/2|+1/2$ and $\tilde{l}_i =|l_i +1/2 |-1/2$ $(i=1,2,3)$.
For  $l_0 ,l_1,l_2, l_3 \in \Zint$, we have polynomials $Q(\lambda ,\mu)$, $P_1 (\lambda ,\mu)$, \dots , $P_6 (\lambda ,\mu)$ such that if $P_2 (\lambda ,\mu)\neq 0$ then there exists a solution $f_{HK}(x ;l_0, l_1, l_2, l_3 ;\lambda ,\mu )$ to Eq.(\ref{eq:Hg}) of the form 
\begin{align}
& f _{HK} (x ;l_0, l_1, l_2, l_3;\lambda ,\mu ) = \exp \left( \kappa x \right) \left( \sum _{i=0}^3 \sum_{j=0}^{\tilde{l}_i-1} \tilde{b} ^{(i)}_j \left( \frac{d}{dx} \right) ^{j} \Phi _i(x, \alpha ) \right)
\label{Lalpha}
\end{align}
for some values $\alpha $, $\kappa$ and $\tilde{b} ^{(i)}_j$ $(i=0,1,2 ,3, \: j= 0,\dots ,\tilde{l}_i-1)$, and the values $\alpha $ and $\kappa $ are expressed as
\begin{align}
& \wp (\alpha )= \frac{P_1 (\lambda ,\mu)}{P_2 (\lambda ,\mu)}, \quad \wp '(\alpha )= \frac{P_3 (\lambda  ,\mu)}{P_4 (\lambda ,\mu)}\sqrt{-Q(\lambda ,\mu)} , \label{eq:alphakappa} \\
& \kappa = \frac{P_5 (\lambda ,\mu)}{P_6 (\lambda ,\mu)}\sqrt{-Q(\lambda ,\mu)} . \nonumber
\end{align}
Regarding the periodicity of the function $f_{HK} (x ) =f_{HK} (x ;l_0, l_1, l_2, l_3;\lambda ,\mu ), $
we have
\begin{align}
& f_{HK} (x+2\omega _j ) = \exp (-2\eta _j \alpha +2\omega _j \zeta (\alpha ) +2 \kappa \omega _j ) f _{HK} (x ), \label{ellint}
\end{align}
for $j=1,3$, where $\eta _j =\zeta (\omega _j)$.
\end{prop}
It follows from Proposition \ref{thm:alpha} that a solution to the Fuchsian differential system $D_{y_1 }( l_1 +1/2, l_2 +1/2, l_3 +1/2, -l_0+1/2; \lambda ,\mu ;k)$ is expressed in the form of the Hermite-Krichever Ansatz for the case $l_0, l_1, l_2 , l_3 \in \Zint $ by setting $z=(\wp (x)-e_1)/(e_2-e_1)$, $t=(e_3-e_2)/(e_1-e_2)$, $y_1 (z)= z^{l_1/2} (z-1)^{l_2/2} (z-t)^{l_3/2}f_{HK}(x;l_0, l_1, l_2, l_3;\lambda ,\mu )$, $y_2 (z)= (dy_1(z)/dz -a_{11}(z)y_1(z))/a_{12}(z)$.

We now consider the Hermite-Krichever Ansatz for the case $l_0=l_1=l_2=l_3=0$ in detail, which was demonstrated in \cite{TakP}.
The differential equation (\ref{eq:Hg}) is written as 
\begin{equation}
\left\{ -\frac{d^2}{dx^2} + \frac{\wp ' (x)}{\wp (x) -\wp (\delta )} \frac{d}{dx} - \frac{ 4\mu \lambda (\lambda -1) (\lambda -t)(e_2-e_1)^2}{\wp (x) -\wp (\delta )} +C \right\} f (x)=0,
\label{Hgkr1l00}
\end{equation}
We assume that $\delta  \not \equiv 0$ mod $\omega _1 \Zint \oplus \omega _3 \Zint$. The condition that the regular singularities $x= \pm \delta $ are apparent is written as
\begin{align}
& C= 2(2 \lambda (\lambda -1) (\lambda -t) \mu^2 -(3\lambda ^2-2(1+t) \lambda +t ) \mu )(e_2-e_1),  \label{pgkrapl00} 
\end{align}
(see Eq.(\ref{pgkrap})).
We consider Eq.(\ref{Hgkr1l00}) with the condition in Eq.(\ref{pgkrapl00}). The polynomial $Q(\lambda ,\mu) $ in Eq.(\ref{eq:alphakappa}) is calculated as
\begin{align}
& Q(\lambda ,\mu)= -2\mu  (2\lambda \mu -1) (2(\lambda -1) \mu -1) (2(\lambda -t) \mu -1)/(e_2-e_1). 
\end{align}
There exists a solution $f_{HK}(x ) (=f_{HK}(x ;0,0,0,0;\lambda ,\mu ))$ to Eq.(\ref{Hgkr1l00}) that can be expressed in the form of the Hermite-Krichever Ansatz as
\begin{align}
& f _{HK} (x) = \bar{b} ^{(0)}_0 \exp (\kappa x) \Phi _0 (x, \alpha ),
\label{eq:HK0000}
\end{align}
if $\mu \neq 0$.
The values $\alpha $ and $\kappa $ are determined as
\begin{align}
& \wp (\alpha )= e_1+(e_2 -e_1) \left(\lambda  - \frac{1}{2 \mu} \right) , \quad \wp '(\alpha )= -\frac{(e_2-e_1)^2\sqrt{-Q(\lambda ,\mu)}}{2\mu ^2} , \\
& \kappa = \frac{(e_2-e_1) \sqrt{-Q(\lambda ,\mu)}}{2\mu } ,\nonumber
\end{align}
and we have
\begin{align}
& \lambda  = \frac{1}{e_2-e_1} \left\{ \wp (\alpha ) -e_1-\frac{\wp ' (\alpha )}{2\kappa } \right\} ,\quad  \mu  = -\frac{(e_2-e_1 )\kappa  }{\wp ' (\alpha )} . \label{eq:lmalka0000}
\end{align}

\section{Integral representation of solutions to Fuchsian system} \label{sec:intrepFs}

We show that solutions to the Fuchsian system (Eq.(\ref{eq:dYdzAzY})) have integral representations for the case $\theta _0, \theta _1, \theta _t, \theta _{\infty } \in \Zint $, $\theta _0 + \theta _1 + \theta _t + \theta _{\infty } \in 1+2\Zint $ by use of the function in the form of Hermite-Krichever Ansatz.
\begin{thm} \label{thm:xinterep}
Assume that $l_0 , l_1 , l_2 ,l_3 \in \Zint$ and $l_0 +l_1 +l_2 +l_3 \in 2\Zint$.
Let $f_{HK} (x)(=f_{HK} (x; l_0 , l_1 , l_2 ,l_3 ;\lambda ,\mu ))$ be the solution expressed in the form of the Hermite-Krichever Ansatz in Proposition \ref{thm:alpha}.
Set 
\begin{align}
& \tilde{\wp }(x)=\frac{\wp (x) -e_1}{e_2-e_1}, \quad \kappa _1=-\frac{l_0+l_1+l_2+l_3+1}{2}.
\end{align}
Then the function $\tilde{y}^{(i)}_1 (z)$ $(i=0,1,2,3)$ defined by
\begin{align}
 \tilde{y}^{(i)} _1 (z)=& \int _{-\tilde{\wp }^{-1} (z )+2\omega _i}^{\tilde{\wp }^{-1}(z)} \left[ \left\{ \frac{\kappa _1 }{e_2- e_1}+ \left( \tilde{\wp }(\xi )- \lambda -\frac{\kappa _1}{\mu }\right) \sum _{j=1}^3 \frac{l_j}{2(\wp (\xi )-e_j)} \right\} \wp ' (\xi ) f _{HK} (\xi ) \right.  \label{eq:intrepy1x}\\
& \quad \quad \quad \left. + \left( \tilde{\wp }(\xi )- \lambda -\frac{\kappa _1}{\mu }\right) \frac{df _{HK} (\xi )}{d\xi } \right] \left( \prod _{j=1}^3 (\wp(\xi )-e_j) ^{l_j/2} \right) \frac{(z-\tilde{\wp }(\xi ))^{\kappa _1 }}{(\tilde{\wp }(\xi )-\lambda )} d\xi , \nonumber 
\end{align}
is a solution to the Fuchsian differential equation $D_{y_1 }(\tilde{\theta }_0, \tilde{\theta }_1, \tilde{\theta }_t, \tilde{\theta }_{\infty}; \tilde{\lambda },\tilde{\mu })$ (see Eq.(\ref{eq:linP6})), where
\begin{align}
& \tilde{\theta }_0=\frac{-l_0 +l _1 -l_2 -l_3}{2}, \quad \tilde{\theta }_1=\frac{-l_0 -l _1 +l_2 -l_3}{2} , \quad \tilde{\theta }_t=\frac{-l_0 -l _1 -l_2 +l_3}{2} ,\\
& \tilde{\theta }_{\infty}=\frac{-l_0 +l _1 +l_2 +l_3}{2}+1 , \quad \tilde{\lambda }= \lambda +\frac{\kappa _1}{\mu}, \quad \tilde{\mu }=\mu . \nonumber
\end{align}
For the case $\kappa _1 \leq -1$, we interpret the integral as a half of the value integrated over the cycle from a point sufficiently close to $\xi = -\tilde{\wp }^{-1} (z )+2\omega _i$, turning around the point $\xi = -\tilde{\wp }^{-1} (z )+2\omega _i$ clockwise, moving to the point sufficiently close to $\xi = \tilde{\wp }^{-1} (z )$, turning around the point $\xi = \tilde{\wp }^{-1} (z )$ anticlockwise and returning to the initial point.
\end{thm}
\begin{proof}
By changing the variable $w=(\wp (\xi) -e_1)/(e_2-e_1)$, substituting $y_1(w)=w^{l_1/2}(w-1)^{l_2/2}(w-t)^{l_3/2} f_{HK}(\xi)$ and multiplying Eq.(\ref{eq:yt1zintrep}) by $t\kappa _2 (u_0+\theta _0)/(-\lambda (\lambda -t)\mu) \tilde{y}_1 (z)$, we obtain the integrand.
We consider the integral contour $[\alpha _{z} ,\alpha _i]$.
Let $o\in  \Cplx \setminus \{0,1,t\}$ be the initial point of the contour in the $w$-plane, and $\pm x_0$ (resp. $\pm x$) be the point such that $\wp (\pm x_0)= o$ (resp. $\wp (\pm x)= z$).
We choose $x_0$ sufficiently close to $x$.
The contour $\alpha _z$ in the $z$-plane corresponds to the contour whose initial point is $x_0$ and turning $x$ anticlockwise and returning either to $x_0$ or to the contour whose initial point is $-x_0$ and turning $-x$ anticlockwise and returning to $-x_0$, depending on the choice of branching.
The contour $\alpha _{\infty }$ (resp. $\alpha _0$, $\alpha _1$, $\alpha _t$) in the $w$-plane corresponds either to the contour whose initial point is $x_0$ and ends at $-x_0$ (resp. $-x_0+2\omega _1$, $-x_0+2\omega _2$, $-x_0+2\omega _3$) or to the reverse contour.
By analytic continuation along the cycle $\alpha _z$, the integrand is multiplied by $-1$ because of the factor $(z-\tilde{\wp }(\xi ))^{\kappa _1 }$ $(\kappa _1 \in \Zint +1/2)$, and the integral tends to zero in the limit as $x_0 \rightarrow x$ for the case $\kappa >-1$.
Hence the contour  $[\alpha _{z} ,\alpha _{\infty }]$ (resp. $[\alpha _{z} ,\alpha _0]$, $[\alpha _{z} ,\alpha _{1}]$, $[\alpha _{z} ,\alpha _{t}]$) corresponds to a contour that runs twice from $x $ to $-x$ (resp. $-x+2\omega _1$, $-x+2\omega _2$, $-x+2\omega _3$).
We therefore obtain the theorem.
\end{proof}
Note that if $l_0, l_1 ,l_2 ,l_3 \in \Zint $ and $l_0 + l_1 +l_2 +l_3 \in 2\Zint +1$, then the integral is equal to zero.
It follows from the assumption $l_0, l_1 ,l_2 ,l_3 \in \Zint $ and $l_0 + l_1 +l_2 +l_3 \in 2\Zint $ that $\tilde{\theta }_0,\tilde{\theta }_1,\tilde{\theta }_t ,\tilde{\theta }_{\infty} \in \Zint $ and $\tilde{\theta }_0+\tilde{\theta }_1+\tilde{\theta }_t +\tilde{\theta }_{\infty} \in 1+2\Zint$.

For given $\theta _0, \theta _1, \theta _t, \theta _{\infty }$ such that $\theta _0, \theta _1, \theta _t, \theta _{\infty } \in \Zint $ and $\theta _0 + \theta _1 + \theta _t + \theta _{\infty } \in 1+2\Zint $,
we have integral representations of solutions by choosing $l_0 , l_1 ,l_2 ,l_3 $ appropriately.
More precisely, we have the following corollary:
\begin{cor} \label{cor:xinterep}
Assume that $\theta _0, \theta _1, \theta _t, \theta _{\infty } \in \Zint $ and $\theta _0 + \theta _1 + \theta _t + \theta _{\infty } \in 1+2\Zint $.
Set 
\begin{align}
& l_0=\frac{-\theta _0- \theta _1- \theta _t- \theta _{\infty }+1}{2},\quad l_1=\frac{\theta _0- \theta _1- \theta _t+ \theta _{\infty }-1}{2}, \\
& l_2=\frac{-\theta _0+ \theta _1- \theta _t+ \theta _{\infty }-1}{2},\quad l_3=\frac{-\theta _0- \theta _1+ \theta _t+ \theta _{\infty }-1}{2}, \nonumber \\
& \tilde{\wp }(x)=\frac{\wp (x) -e_1}{e_2-e_1}, \quad \tilde{\kappa }_1=\frac{\theta _0+ \theta _1+ \theta _t- \theta _{\infty }}{2}. \nonumber
\end{align}
Let $f_{HK} (x)=(f_{HK} (x; l_0 , l_1 , l_2 ,l_3 ;\lambda -\tilde{\kappa _1}/\mu, \mu))$ be the solution expressed in the form of the Hermite-Krichever Ansatz in Proposition \ref{thm:alpha}.
Then the function $\tilde{y}^{(i)}_1 (z)$ $(i=0,1,2,3)$ defined by
\begin{align}
\tilde{y}^{(i)} _1 (z)= &\int _{-\tilde{\wp }^{-1} (z )+2\omega _i}^{\tilde{\wp }^{-1}(z)} \left[ \left\{ \frac{\tilde{\kappa _1 }}{e_2- e_1}+ \left( \tilde{\wp }(\xi )- \lambda \right) \sum _{j=1}^3 \frac{l_j}{2(\wp (\xi )-e_j)} \right\} \wp ' (\xi ) f _{HK} (\xi ) \right. \label{eq:y1tcor} \\
& \quad \quad \left. + \left( \tilde{\wp }(\xi )- \lambda \right) \frac{df _{HK} (\xi )}{d\xi } \right] \prod _{j=1}^3 (\wp(\xi )-e_j) ^{l_j/2}  \frac{(z-\tilde{\wp }(\xi ))^{\tilde{\kappa }_1 }}{(\tilde{\wp }(\xi )-\lambda +\tilde{\kappa _1}/\mu)} d\xi ,\nonumber
\end{align}
is a solution to the Fuchsian differential equation $D_{y_1 }(\theta _0, \theta _1, \theta _t, \theta _{\infty}; \lambda ,\mu )$.
For the case $\kappa _1 \leq -1$, we interpret the integral as the one in Theorem \ref{thm:xinterep}.
\end{cor}
Let $a_{11}(z), a_{12}(z)$ be the functions defined in Eq.(\ref{eq:dYdzAzY}) and set $\tilde{y}_2 ^{(i) }(z) = ( d{\tilde y}_1^{(i) }(z)/dz -a_{11}(z)\tilde{y}_1^{(i) }(z))/a_{12}(z)$ $(i=0,1,2,3)$.
Then the function $Y={}^t ( \tilde{y}_1 ^{(i) }(z) ,\tilde{y}_2 ^{(i) }(z) )$ is a solution to the Fuchsian differential system $D_{Y}(\theta _0, \theta _1, \theta _t, \theta _{\infty}; \lambda ,\mu ;k)$ (see Eq.(\ref{eq:dYdzAzY})).
Note that the function $\tilde{y}_2 ^{(i) }(z)$ is also expressed as the form like Eq.(\ref{eq:y1tcor}) by combining the expression in Eq.(\ref{eq:yt1zintrep}), the relation $y_1 (w)= w^{l_1/2} (w-1)^{l_2/2} (w-t)^{l_3/2}f_{HK}(\tilde{\wp }^{-1}(w ))$ and the relation among $dy_2 (w)/dw$, $dy_1(w)/dw$ and $y_1(w)$.

We can calculate the monodromy of the Fuchsian system $D_{Y}(\theta _0, \theta _1, \theta _t, \theta _{\infty}; \lambda ,\mu ;k)$ and the Fuchsian equation $D_{y_1 }(\theta _0, \theta _1, \theta _t, \theta _{\infty}; \lambda ,\mu )$ for the case $\theta _0, \theta _1, \theta _t, \theta _{\infty } \in \Zint $ and $\theta _0 + \theta _1 + \theta _t + \theta _{\infty } \in 1+2\Zint $ in principal by considering the integral representations of solutions and their asymptotics around the singularities.
We will do this for the case $(\theta _0, \theta _1, \theta _t, 1-\theta _{\infty }) =(0,0,0,0)$ in the next section.
\section{Integral representation of solutions to the Fuchsian equation for the case $(\theta _0, \theta _1, \theta _t, 1-\theta _{\infty }) =(0,0,0,0)$} \label{sec:intrepFs0000}

We consider the integral representation of solutions to the Fuchsian equation for the case $(\theta _0, \theta _1, \theta _t, 1-\theta _{\infty }) =(0,0,0,0)$.
For this case, the function $f_{HK}(\xi )$ in the integrand of Eq.(\ref{eq:intrepy1x}) is written in the form of the Hermite-Krichever Ansatz for the case $l_0=l_1=l_2=l_3=0$, and it is described by Eq.(\ref{eq:HK0000}).
The values $\lambda ,\mu $ for the case $l_0=l_1=l_2=l_3=0$ and the values $\alpha , \kappa $ are related by Eq.(\ref{eq:lmalka0000}).
By substituting Eq.(\ref{eq:lmalka0000}) into Eq.(\ref{eq:y1tcor}) and the integral representation of the function $\tilde{y}_2 ^{(i) }(z)$ like Eq.(\ref{eq:y1tcor}), multiplying by appropriate constants and applying the formula $\wp (x) -\wp (\xi )=-\sigma (x+ \xi) \sigma (x- \xi)/(\sigma (x)^2 \sigma (\xi)^2)$,
we have the following proposition:
\begin{prop} \label{prop:intrep0001}
Set
\begin{align}
& f_i(x)=\int_{-x+2\omega _i}^x \frac{e^{(\kappa +\zeta (\alpha ))\xi }\sigma (x)\sigma (\xi -\alpha )}{\sqrt{\sigma (x-\xi )\sigma (x+\xi )}}d\xi , \label{eq:intrep0001} \\
& g_i(x)= \frac{1}{4k (e_2-e_1)} \int_{-x+2\omega _i}^x \left( \kappa ^2 +\frac{\wp '(\xi ) +\wp '(\alpha )}{\wp (\xi ) -\wp (\alpha )} \kappa +2\wp (\xi ) +\wp (\alpha ) \right) \nonumber \\
& \quad \quad \quad \quad \quad \quad \quad \quad \quad \quad \quad \quad \quad \quad \quad \quad \quad \quad \frac{e^{(\kappa +\zeta (\alpha ))\xi }\sigma (x)\sigma (\xi -\alpha )}{\sqrt{\sigma (x-\xi )\sigma (x+\xi )}}d\xi . \nonumber
\end{align}
The function ${}^t(f_i (x) , g_i (x))$ $(i=0,1,2,3, \; z=(\wp (x)-e_1)/(e_2-e_1))$ is a solution to the Fuchsian differential system $D_{Y} (0,0,0,1;\lambda ,\mu ; k)$, where 
\begin{align}
& \lambda =\frac{\wp (\alpha ) -e_1}{e_2-e_1}, \quad \mu= -\frac{(e_2-e_1 )\kappa }{\wp ' (\alpha )} .
\end{align}
In particular, the function $f_i (x)$ $(i=0,1,2,3, \; z=(\wp (x)-e_1)/(e_2-e_1))$ is a solution to the Fuchsian differential equation $D_{y_1} (0,0,0,1;\lambda ,\mu)$ and the differential equation can also be written as
\begin{align}
& \frac{d^2y}{dx^2} +\left\{ \left( \sum _{j=1}^3 \frac{1}{2(\wp (x)- e_j)} \right) -\frac{1}{\wp (x) -\wp (\alpha )} \right\} \wp '(x) \frac{dy}{dx} \\
& \quad \quad \quad +\left\{ -\kappa ^2-\frac{\wp ' (\alpha )}{\wp (x) -\wp (\alpha )}\kappa +\wp (x)-\wp (\alpha ) \right\} y=0.\nonumber
\end{align}
\end{prop}
The monodromy matrix for the Fuchsian differential system $D_{Y} (0,0,0,1;\lambda ,\mu ; k)$ with respect to a basis of solutions 
$ \left\{ \left( \begin{array}{cc}
y^{\{1 \} }_1(z) \\
y^{\{1 \} }_2(z)
\end{array}
\right)  ,
\left( \begin{array}{cc}
y^{\{2 \} }_1 (z) \\
y^{\{2 \} }_2 (z)
\end{array}
\right) \right\} $
along a cycle $\gamma $ coincides with the monodromy matrix for the Fuchsian differential equation $D_{y_1} (0,0,0,1;\lambda ,\mu )$ with respect to a basis of solutions 
$ \{ y^{\{1 \} }_1(z) ,y^{\{ 2 \} }_1(z)\} $
along the cycle $\gamma $.
Hence we investigate the monodromy matries for the Fuchsian differential equation $D_{y_1} (0,0,0,1;\lambda ,\mu)$ by applying integral representations of solutions $f_i(x)$ $(i=0,1,2,3)$.

Assume that $\alpha \not \equiv 0$ mod $\omega _1\Zint \oplus \omega _3 \Zint$.
By considering the exponents of the singularities, we have a basis of local solutions to the Fuchsian differential equation $D_{y_1} (0,0,0,1;\lambda ,\mu)$ about $x=0$ and $\omega _i$ $(i=1,2,3)$ of the form
\begin{align}
& s^{(0)}_1 (x)= x+c^{(0)}_2 x^2+ \dots ,  \quad s^{(0)}_2 (x)= s^{(0)}_1 (x) \log x +\tilde{c}^{(0)}_1 x + \tilde{c}^{(0)}_2  x^2 + \dots ,\\
& s^{(i)}_1 (x)= 1+c^{(i)}_1 (x-\omega _i) + \dots ,\quad s^{(i)}_2 (x)= s^{(i)}_1 (x) \log (x-\omega _i)  +\tilde{c}^{(i)}_0 + \tilde{c}^{(i)}_1 (x -\omega _i) + \dots . \nonumber
\end{align}
Let $\gamma _i$ $(i=0,1,2,3)$ be the cycle turning anti-clockwise around $x=\omega _i$, and $f^{\gamma }(x)$ be the function which is continued analytically along the cycle $\gamma $.
Then we have
\begin{align}
& (s^{(i),\gamma _i}_1 (x), s^{(i),\gamma _i}_2 (x)) =(s^{(i)}_1 (x), s^{(i)}_2 (x)) \left(
\begin{array}{cc}
1 &  2\pi \sqrt{-1} \\
0 & 1 
\end{array}
\right)
\quad (i=0,1,2,3). \label{eq:silocmonod}
\end{align}
We now relate the local solutions $(s^{(i)}_1 (x), s^{(i)}_2 (x))$ $(i=0,1,2,3)$ to the solutions of integral representations $f_0 (x), \dots ,f_3(x)$.
Since $f_0(x)$ is a solution to the Fuchsian differential equation $D_{y_1} (0,0,0,1,\lambda ,\mu)$, it is expressed as a linear combination of $s^{(0)}_1 (x)$ and $s^{(0)}_2 (x)$.
We set $\xi =x \nu $. Since $\lim _{x \rightarrow 0} \sigma (x)/x =1$, we have the following asymptotic limit as $x \rightarrow 0$:
\begin{align}
& f_0(x)= \int _{-1} ^1 \frac{e^{x(\kappa +\zeta (\alpha ))\nu }\sigma (x)\sigma (x \nu -\alpha )}{\sqrt{\sigma (x(1+ \nu ))\sigma (x(1- \nu ))}} xd\nu  \sim \int _{-1} ^1 \frac{x^2 \sigma (-\alpha )}{x\sqrt{(1+\nu )(1-\nu )}}d\nu  = \sigma (-\alpha )\pi x.
\end{align}
Hence we have
\begin{equation}
f_0(x)= \sigma (-\alpha ) \pi s^{(0)}_1 (x).
\label{eq:f0s0}
\end{equation}
We consider the asymptotics of $f_0(x)$ in the limit as $x \rightarrow \omega _i$ $(i=1,2,3)$.
By using the formula $\sigma (x+2\omega _i)= -\sigma (x)\exp (2\eta _i (x+\omega _i))$, we have
\begin{align}
& f_0(x)=\int _0^1 
\frac{e^{x(\kappa +\zeta (\alpha ))\nu } \sigma (x) \sigma (x \nu  -\alpha )xd\nu }{\sqrt{-e^{2\eta _i(x(1+ \nu )-\omega _i)}\sigma (x(1+ \nu )-2\omega _i)\sigma (x(1- \nu ))}} \\
& \quad +\int _{-1}^0 \frac{e^{x(\kappa +\zeta (\alpha ))\nu }\sigma (x)\sigma (x \nu -\alpha )xd\nu }{\sqrt{-e^{2\eta _i(x(1- \nu )-\omega _i)}\sigma (x(1- \nu )-2\omega _i)\sigma (x(1+ \nu ))}} \nonumber \\
& \quad \sim \int _0^1 \frac{e^{x(\kappa +\zeta (\alpha ))-\eta _i(2x-\omega _i)}\sigma (x)\sigma (x -\alpha )xd\nu }{x\sqrt{-(1+ \nu -2\omega _i/x ) (1- \nu )}} +\int _{-1}^0 \frac{e^{-x(\kappa +\zeta (\alpha ))-\eta _i(2x-\omega _i)}\sigma (x)\sigma (-x -\alpha )xd\nu }{x\sqrt{-(1- \nu -2\omega _i/x ) (1+ \nu )}} \nonumber \\
& \quad \sim -\log (\omega _i-x) e^{\omega _i (\kappa +\zeta (\alpha )-\eta _i )} \sigma (\omega _i) \sigma (\omega _i- \alpha ) -\log (\omega _i-x) e^{\omega _i (-(\kappa +\zeta (\alpha ))-\eta _i )} \sigma (\omega _i) \sigma (-\omega _i- \alpha ) \nonumber \\
& \quad =  -\log (\omega _i-x) e^{\omega _i (\kappa +\zeta (\alpha )- \eta _i) } \sigma (\omega _i) \sigma (\omega _i- \alpha ) (1 - e^{-2\omega _i (\kappa +\zeta (\alpha ))+2\eta _i \alpha }) .\nonumber
\end{align}
Since $f_0(x)$ is a solution to Eq.(\ref{eq:dYdzAzY}), it can be expressed as a linear combination of $s^{(i)}_1 (x)$ and $s^{(i)}_2 (x)$, and we have
\begin{equation}
f_0(x)= -e^{\omega _i (\kappa +\zeta (\alpha )- \eta _i )} \sigma (\omega _i) \sigma (\omega _i- \alpha ) (1- e^{-2\omega _i (\kappa +\zeta (\alpha ))+2\eta _i \alpha }) s^{(i)}_2 (x) +c^{(0,i)} s^{(i)}_1 (x),
\label{eq:f0si}
\end{equation}
for some constant $c^{(0,i)} $.
Next, we express the function $f_i(x)$ $(i=1,2,3)$ as a linear combination of $s^{(j)}_1 (x)$ and $s^{(j)}_2 (x)$ for $j\in \{ 0,1,2,3\}$.
We set $\xi =(\omega _i-x) \nu  +\omega _i$, whereupon we have
\begin{align}
& f_i(x)= \int _{-1} ^1 \frac{e^{x(\kappa +\zeta (\alpha ))((\omega _i-x) \nu  +\omega _i) }\sigma (x)\sigma ((\omega _i-x) \nu  +\omega _i-\alpha )}{\sqrt{\sigma ((\omega _i-x) (1+\nu ))\sigma ((\omega _i-x) (1-\nu )-2\omega _i )}} (x-\omega _i)d\nu .
\end{align}
Similarly, we have
\begin{align}
& f_i(x) \sim \sqrt{-1} \pi \sigma (\omega _i ) \sigma (\omega _i -\alpha ) e^{\omega _i(\kappa +\zeta (\alpha )-\eta _i)} , \quad (x \rightarrow \omega _i), \\
& f_i(x) \sim - \sqrt{-1} \sigma ( -\alpha ) (1-e^{2\omega _i (\kappa +\zeta (\alpha ))-2\eta _i \alpha }) x \log x , \quad (x \rightarrow 0) , \nonumber
\end{align}
and
\begin{align}
& f_i(x) \sim \sigma (\omega _j ) \sigma (\omega _j- \alpha ) e^{\omega _j (\kappa +\zeta (\alpha )-\eta _j ) }(1-e^{2(\omega _i- \omega _j )(\kappa +\zeta (\alpha ))+2(\eta _j -\eta _i) \alpha }) \log (\omega _j -x) ,
\end{align}
in the limit as $x \rightarrow \omega _j $,  $(j\neq 0,i)$, where we have used  Legendre's relation, $\eta_i \omega _j-\eta _j \omega _i=\pm \pi \sqrt{-1}/2$.
Therefore we have
\begin{align}
 f_i(x) & = \sqrt{-1} \pi \sigma (\omega _i ) \sigma (\omega _i -\alpha ) e^{\omega _i(\kappa +\zeta (\alpha )-\eta _i)} s^{(i)}_1 (x) \label{eq:fisi} \\
& =-\sqrt{-1} \sigma ( -\alpha ) (1-e^{2\omega _i (\kappa +\zeta (\alpha ))-2\eta _i \alpha })s^{(0)}_2 (x) +c^{(i,0)} s^{(0)}_1 (x) \nonumber \\
& =\sigma (\omega _j ) \sigma (\omega _j- \alpha ) e^{\omega _j (\kappa +\zeta (\alpha )-\eta _j ) }(1-e^{2(\omega _i- \omega _j )(\kappa +\zeta (\alpha ))+2(\eta _j -\eta _i) \alpha }) s^{(j)}_2 (x) +c^{(i,j)} s^{(j)}_1 (x), \nonumber
\end{align}
for some constants $c^{(i,0)} $ and $c^{(i,j)} $.

We consider the monodromy matrices on the basis $(f_0(x) , f_1(x))$. Set
\begin{equation}
e[i]=\exp( 2\omega _i (\kappa +\zeta (\alpha ))-2\eta _i \alpha ), \quad  (i=1,2,3).
\end{equation}
It follows from Eqs.(\ref{eq:silocmonod}, \ref{eq:f0s0}, \ref{eq:fisi}) that
\begin{align}
& (f^{\gamma _0}_0(x) , f^{\gamma _0}_1(x))= (\sigma (-\alpha ) \pi s^{(0),\gamma _0}_1 (x), -\sqrt{-1} \sigma ( -\alpha ) (1-e[1])s^{(0),\gamma _0}_2 (x) +c^{(i,0)} s^{(0),\gamma _0}_1 (x)) \label{eq:f0f1gam0} \\
& = (\sigma (-\alpha ) \pi s^{(0)}_1(x), -\sqrt{-1} \sigma ( -\alpha ) (1-e[1])(s^{(0)}_2 (x) + 2\pi \sqrt{-1} s^{(0)}_1 (x) )+c^{(i,0)} s^{(0)}_1 (x)) \nonumber \\
& = (f_0(x) , f_1(x))
\left(
\begin{array}{cc}
1 &  2(1-e[1]) \\
0 & 1 
\end{array}
\right) . \nonumber 
\end{align}
Similarly, it follows from Eqs.(\ref{eq:silocmonod}, \ref{eq:f0si}, \ref{eq:fisi}) that
\begin{align}
& (f^{\gamma _1}_0(x) , f^{\gamma _1}_1(x))=  (f_0(x) , f_1(x))
\left(
\begin{array}{cc}
1 &  0 \\
-2(1-1/e[1]) & 1 
\end{array}
\right) . \label{eq:f0f1gam1}
\end{align}
If $e[1] \neq 0$, then it follows from the asymptotic limits as $x \rightarrow 0$ and $x\rightarrow \omega _1$ that the functions $f_0(x)$ and $f_1(x)$ form a basis of solutions to Eq.(\ref{eq:dYdzAzY}),
and the functions $f_j(x)$ $(j=2,3)$ are written as linear combinations of $f_0(x)$ and $f_1(x)$.
Write $f_j(x)= \tilde{c}_{0,j}f_0(x) +\tilde{c}_{1,j}f_1(x)$. Then the coefficients $\tilde{c}_{0,j}, \tilde{c}_{1,j}$ are determined by considering the asymptotic limits as $x \rightarrow \omega _1$ and $x\rightarrow 0$, and we have
\begin{align}
& \tilde{c}_{0,j} =  \frac{e[1]-e[j]}{1-e[1]} , \quad  \tilde{c}_{1,j} = \frac{1-e[j]}{1-e[1]} .
\end{align}
Therefore
\begin{align}
& (f^{\gamma _j}_0(x) , f^{\gamma _j}_1(x))=  (f^{\gamma _j}_0(x) , f^{\gamma _j}_j(x))
\left(
\begin{array}{cc}
1 & -\tilde{c}_{0,j}/\tilde{c}_{1,j} \\
0 & 1/\tilde{c}_{1,j} 
\end{array}
\right) \label{eq:f0f1gamj} \\
& = (f_0(x) , f_j(x))
\left(
\begin{array}{cc}
1 &  0 \\
-2(1-1/e[j]) & 1 
\end{array}
\right)  
\left(
\begin{array}{cc}
1 & -\tilde{c}_{0,j}/\tilde{c}_{1,j} \\
0 & 1/\tilde{c}_{1,j} 
\end{array}
\right) \nonumber \\
& =(f_0(x) , f_1(x))
\left(
\begin{array}{cc}
1 & \tilde{c}_{0,j} \\
0 & \tilde{c}_{1,j} 
\end{array}
\right) 
\left(
\begin{array}{cc}
1 &  0 \\
-2(1-1/e[j]) & 1 
\end{array}
\right)
\left(
\begin{array}{cc}
1 & -\tilde{c}_{0,j}/\tilde{c}_{1,j} \\
0 & 1/\tilde{c}_{1,j} 
\end{array}
\right) \nonumber \\
&=  (f_0(x) , f_1(x))
\left(
\begin{array}{cc}
1+2\frac{(e[1]-e[j])(e[j]-1)}{(e[1]-1)e[j]} &  2\frac{(e[1]-e[j])^2}{(e[1]-1)e[j]} \\
-2\frac{(e[j]-1)^2}{(e[1]-1)e[j]} & 1-2\frac{(e[1]-e[j])(e[j]-1)}{(e[1]-1)e[j]} 
\end{array}
\right) ,\nonumber 
\end{align}
and we have obtained the monodromy matrices for the basis $(f_0(x), f_1(x))$ on the cycles $\gamma _2$, $\gamma _3 $.

We consider the monodromy preserving deformation with respect to the basis $(f_0(x), f_1(x))$.
Assume that the values $e[1]$, $e[3]$ are preserved while varying the ratio $\omega _3/\omega _1$.
Then the monodromy is preserved by Eqs.(\ref{eq:f0f1gam0}, \ref{eq:f0f1gam1}, \ref{eq:f0f1gamj}) and the equality $e[1]+e[2]+e[3]=0$.
Since the values $e[1]$, $e[3]$ are preserved by monodromy preserving deformation, we have 
\begin{align}
& -2\eta _1 \alpha +2\omega _1 \zeta (\alpha ) +2 \kappa \omega _1 = \pi \sqrt{-1} C_1, \\
& -2\eta _3 \alpha +2\omega _3 \zeta (\alpha ) +2 \kappa \omega _3 = \pi \sqrt{-1} C_3, \nonumber
\end{align}
for constants $C_1$ and $C_3$.
By Legendre's relation, $\eta _1\omega _3-\eta _3\omega _1=\pi \sqrt{-1}/2$, we have 
\begin{align}
& \alpha  = C_3 \omega _1 -C_1 \omega _3  \label{al00},\\
& \kappa = \zeta (C_1 \omega _3 -C_3 \omega _1 ) +C_3 \eta _1 -C_1 \eta _3  , \nonumber
\end{align}
Recall that the sixth Painlev\'e equation has an elliptical representation (see Eq.(\ref{eq:P6ellip})),
and it is a differential equation on $\delta $ with respect to the variable $\tau =\omega_3/\omega _1$.
For the case $(\theta _0, \theta _1, \theta _t, 1-\theta _{\infty }) =(0,0,0,0)$, this equation is written as $d^2\delta /d\tau ^2 =0$.
The variables $\lambda $ and $\delta $ are related by $\lambda =(\wp (\delta )-e_1)/(e_2-e_1)$.
With regards to the integral representation of solutions to $D_{y_1}(0,0,0,1;\lambda ,\mu)$, we have the relations
\begin{align}
& \lambda  = \frac{1}{e_2-e_1} \left\{ \wp (\alpha ) -e_1\right\} ,\quad  \mu  = -\frac{(e_2-e_1 )\kappa  }{\wp ' (\alpha )}  . \label{eq:lmalka0000t}
\end{align}
Hence $\alpha $ plays the role of $\delta $ mod $2\omega _1 \Zint \oplus 2\omega _3 \Zint $, and Eq.(\ref{al00}) corresponds to Picard's solution to the sixth Painlev\'e equation for the case $(\theta _0, \theta _1, \theta _t, 1-\theta _{\infty }) =(0,0,0,0)$ by setting $\omega _1=1/2$ and $\omega _3=\tau /2$.
We therefore reproduce Picard's solution by determining the monodromy of the corresponding Fuchsian equation.

\section{Integral representation of solutions to Heun's equation} \label{sec:intrepHeun}
In section \ref{sec:intrepFs}, we obtained that, if $\theta _0, \theta _1, \theta _t, \theta _{\infty } \in \Zint $ and $\theta _0 + \theta _1 + \theta _t + \theta _{\infty } \in 1+2\Zint $, then we have integral representations of solutions to the Fuchsian equation $D_{y_1}(\theta _0, \theta _1, \theta _t, \theta _{\infty}; \lambda ,\mu )$ (see Eq.(\ref{eq:intrepy1x})) and the Fuchsian system $D_Y(\theta _0, \theta _1, \theta _t, \theta _{\infty}; \lambda ,\mu ,k)$.
In this section we obtain integral representations of solutions to Heun's equation by a suitable choice of the parameters $\lambda $ and $\mu $.

Recall that Heun's differential equation is defined by
\begin{equation}
\frac{d^2y}{dz^2} + \left( \frac{\gamma}{z}+\frac{\delta }{z-1}+\frac{\epsilon}{z-t}\right) \frac{dy}{dz} +\frac{\alpha \beta z -q}{z(z-1)(z-t)} y=0,
\label{eq:Heuns7}
\end{equation}
with the condition $\gamma +\delta +\epsilon =\alpha +\beta +1$.
This equation has an elliptical representation:
Set
\begin{align}
& z=\frac{\wp (x) -e_1}{e_2-e_1}, \quad  t=\frac{e_3-e_1}{e_2-e_1}, \quad f(x)= y z^{\frac{-l_1}{2}}(z-1)^{\frac{-l_2}{2}}(z-t)^{\frac{-l_3}{2}},
\end{align}
then Heun's equation (Eq.(\ref{eq:Heuns7})) is transformed to
\begin{equation}
\left(-\frac{d^2}{dx^2} + \sum_{i=0}^3 l_i(l_i+1)\wp (x+\omega_i)-E\right)f(x)=0,
\label{InoEF0}
\end{equation}
where 
\begin{align}
& l_0= \alpha -\beta -1/2,\quad l_1= -\gamma +1/2, \quad l_2=-\delta +1/2, \quad l_3=-\epsilon +1/2, \\
& E=(e_2-e_1)(-4q+(-(\alpha -\beta)^2 +2\gamma ^2+6\gamma \epsilon +2\epsilon ^2 -4\gamma -4\epsilon  -\delta ^2 +2\delta +1)/3 \nonumber \\
& \quad \quad +(-(\alpha -\beta ) ^2 +2\gamma ^2+6\gamma \delta +2\delta^2 -4\gamma -4\delta -\epsilon ^2+2\epsilon +1)t/3). \nonumber
\end{align}

We obtained in section \ref{sec:intrepFs} that, if $l_0 , l_1 , l_2 ,l_3 \in \Zint$ and $l_0 +l_1 +l_2 +l_3 \in 2\Zint$, then the function $\tilde{y}^{(i)}_1 (z) $
defined by 
\begin{align}
 \tilde{y}^{(i)} _1 (z)=& \int _{-\tilde{\wp }^{-1} (z )+2\omega _i}^{\tilde{\wp }^{-1}(z)} \left[ \left\{ \frac{\kappa _1 }{e_2- e_1}+ \left( \tilde{\wp }(\xi )- \lambda -\frac{\kappa _1}{\mu }\right) \sum _{i=1}^3 \frac{l_i}{2(\wp (\xi )-e_i)} \right\} \wp ' (\xi ) f _{HK} (\xi ) \right.  \label{eq:intrepy1x0}\\
& \quad \quad \quad \left. + \left( \tilde{\wp }(\xi )- \lambda -\frac{\kappa _1}{\mu }\right) \frac{df _{HK} (\xi )}{d\xi } \right] \left( \prod _{i=1}^3 (\wp(\xi )-e_i) ^{l_i/2} \right) \frac{(z-\tilde{\wp }(\xi ))^{\kappa _1 }}{(\tilde{\wp }(\xi )-\lambda )} d\xi , \nonumber 
\end{align}
$(i=0,1,2,3)$ is a solution to the Fuchsian differential equation $D_{y_1 }(\tilde{\theta }_0, \tilde{\theta }_1, \tilde{\theta }_t, \tilde{\theta }_{\infty}; \lambda +\kappa _1 /\mu ,\mu )$, where 
\begin{align}
& \tilde{\theta }_0=\frac{-l_0 +l _1 -l_2 -l_3}{2}, \quad \tilde{\theta }_1=\frac{-l_0 -l _1 +l_2 -l_3}{2} , \quad \tilde{\theta }_t=\frac{-l_0 -l _1 -l_2 +l_3}{2} ,\\
& \tilde{\theta }_{\infty}=\frac{-l_0 +l _1 +l_2 +l_3}{2}+1 , \quad \tilde{\wp }(x)=\frac{\wp (x) -e_1}{e_2-e_1}, \quad \kappa _1=-\frac{l_0+l_1+l_2+l_3+1}{2}, \nonumber
\end{align}
and the function $f_{HK}(x)$ is defined in Theorem \ref{thm:xinterep}.
The Fuchsian equation $D_{y_1 }(\tilde{\theta }_0, \tilde{\theta }_1, \tilde{\theta }_t, \tilde{\theta }_{\infty}; \lambda +\kappa _1 /\mu ,\mu )$ has an apparent singularity at $z=\lambda +\kappa _1/\mu $.
We consider the confluence of the apparent singularity $z =\lambda +\kappa _1/\mu $ to the regular singularity $z=\infty$.
Set $\mu= 0$.
Then the Fuchsian equation is written as Heun's equation
\begin{align}
& \frac{d^2y}{dz^2} + \left( \frac{1-\tilde{\theta }_0}{z}+\frac{1-\tilde{\theta }_1}{z-1}+\frac{1-\tilde{\theta }_t}{z-t}\right) \frac{dy}{dz} \label{Heuninfty} \\
& +\frac{\tilde{\kappa }_1 (\tilde{\kappa }_2+2) z + \tilde{\kappa }_1 (1 -\tilde{\theta }_{\infty} )\lambda  -\tilde{\kappa }_1  ( (\tilde{\kappa }_2 +\tilde{\theta }_t +1)t +(\tilde{\kappa }_2 +\tilde{\theta }_1 +1)) }{z(z-1)(z-t)} y=0 ,\nonumber
\end{align}
where $\tilde{\kappa }_1 =( \tilde{\theta }_{\infty } -\tilde{\theta }_0-\tilde{\theta }_1-\tilde{\theta }_t)/2$ and $\tilde{\kappa }_2 =-( \tilde{\theta }_{\infty } +\tilde{\theta }_0+\tilde{\theta }_1+\tilde{\theta }_t)/2$.
We have $1-\tilde{\theta }_0, 1-\tilde{\theta }_1,1-\tilde{\theta }_t , \tilde{\kappa }_1 +1/2, \tilde{\kappa }_2+2+1/2 \in \Zint$.
For the case $\tilde{\theta }_{\infty} =1$, we set $\mu =bs^2$, $\lambda = c/s$ and consider the limit $s \rightarrow 0$.
Then we have
\begin{align}
& \frac{d^2y}{dz^2} + \left( \frac{1-\tilde{\theta }_0}{z}+\frac{1-\tilde{\theta }_1}{z-1}+\frac{1-\tilde{\theta }_t}{z-t}\right) \frac{dy}{dz} \label{Heuninftytti1} \\
& +\frac{\tilde{\kappa }_1 (\tilde{\kappa }_2+2) z  + \tilde{\kappa }_1 bc^2 - \tilde{\kappa }_1  ( (\tilde{\kappa }_2 +\tilde{\theta }_t +1)t +(\tilde{\kappa }_2 +\tilde{\theta }_1 +1)) }{z(z-1)(z-t)} y=0 ,\nonumber
\end{align}
The following theorem follows from Eq.(\ref{Heuninfty}) by substituting the parameters as indicated:
\begin{thm} \label{thm:interepHeun}
(i) Assume that $\gamma +\delta +\epsilon =\alpha +\beta +1$, $\gamma ,\delta ,\epsilon ,\alpha +1/2 , \beta +1/2 \in \Zint$.
Set 
\begin{align}
& \tilde{l}_0=\alpha -3/2, \; \tilde{l}_1=\delta +\epsilon -\alpha -1/2, \; \tilde{l}_2=\gamma +\epsilon -\alpha -1/2, \; \tilde{l}_3=\gamma + \delta -\alpha -1/2.
\end{align}
Let $f_{HK}(x)= f_{HK} (x; \tilde{l}_0,\tilde{l}_1,\tilde{l}_2,\tilde{l}_3;\lambda ,\mu )$ be the function expressed in the form of the Hermite-Krichever Ansatz.
Set $\tilde{\wp }(x)= (\wp (x)-e_1)/(e_2-e_1)$ and
\begin{align}
& F(\xi; \lambda ,\mu, m)= \mu ^m \left[ \left\{ \frac{-\beta }{e_2- e_1}+ \left( \tilde{\wp }(\xi )- \lambda +\frac{\beta }{\mu }\right) \sum _{i=1}^3 \frac{\tilde{l}_i}{2(\wp (\xi )-e_i)} \right\} \wp ' (\xi ) f _{HK} (\xi ) \right.  \label{eq:intrepy1x00}\\
& \quad \quad \quad \quad \quad \quad \quad \quad \quad \left. + \left( \tilde{\wp }(\xi )- \lambda +\frac{\beta }{\mu }\right) \frac{df _{HK} (\xi )}{d\xi } \right] \left( \prod _{i=1}^3 (\wp(\xi )-e_i) ^{\tilde{l}_i/2} \right) .\nonumber
\end{align}
If $\alpha -\beta \neq 1$ (resp. $\alpha -\beta = 1$) and the integrand in Eq.(\ref{eq:intrepHeun}) (resp. Eq.(\ref{eq:intrepHeun00})) has a non-zero finite limit as $\mu \rightarrow 0$  (resp. $s \rightarrow 0$) for some $m$,
then the functions 
\begin{align}
& \tilde{y}^{(i)}_1 (z)= \int _{-\tilde{\wp }^{-1}(z) +2\omega _i}^{\tilde{\wp }^{-1}(z) } \lim _{\mu \rightarrow 0} F(\xi; \lambda ,\mu , m) \frac{(z-\tilde{\wp }(\xi ))^{-\beta  }}{(\tilde{\wp }(\xi )-\lambda )} d\xi , \label{eq:intrepHeun} \\
& \quad \quad  \lambda  = \frac{t(\alpha -\epsilon )+(\alpha -\delta )-q/\beta }{\alpha -\beta - 1} , \quad (\alpha -\beta \neq 1), \nonumber \\
& \tilde{y}^{(i)}_1 (z)= \int _{-\tilde{\wp }^{-1}(z) +2\omega _i}^{\tilde{\wp }^{-1}(z) } \lim _{s \rightarrow 0} F(\xi; c/s ,b s^2 , m) \frac{(z-\tilde{\wp }(\xi ))^{-\beta  }}{(\tilde{\wp }(\xi )-c/s )} d\xi , \label{eq:intrepHeun00} \\
& \quad \quad bc^2 = t(\alpha -\epsilon )+(\alpha -\delta )+q/(1- \alpha ) , \quad (\alpha -\beta = 1), \nonumber 
\end{align}
$(i=0,1,2,3)$ are solutions to Heun's equation (Eq.(\ref{eq:Heuns7})).\\
(ii) Assume that $l_0,l_1,l_2,l_3 \in \Zint +1/2$ and $l_0+l_1+l_2+l_3 \in 1+2\Zint $.
Set
\begin{align}
& \tilde{l}= \frac{l_0+l_1+l_2+l_3}{2}, \quad \tilde{l}_0=\frac{l_0-l_1-l_2-l_3}{2}-1, \\
& \tilde{l}_1=\frac{-l_0+l_1-l_2-l_3}{2}, \quad \tilde{l}_2=\frac{-l_0-l_1+l_2-l_3}{2}, \quad \tilde{l}_3=\frac{-l_0-l_1-l_2+l_3}{2}, \nonumber \\
& F(\xi; \lambda ,\mu, m)= \mu ^m \left[ \left\{ \frac{\tilde{l}}{e_2- e_1}+ \left( \tilde{\wp }(\xi )- \lambda -\frac{\tilde{l}}{\mu }\right) \sum _{i=1}^3 \frac{\tilde{l}_i}{2(\wp (\xi )-e_i)} \right\} \wp ' (\xi ) f _{HK} (\xi ) \right. \nonumber \\
& \quad \quad \quad \quad \quad \quad \quad \quad \left. + \left( \tilde{\wp }(\xi )- \lambda -\frac{\tilde{l}}{\mu }\right) \frac{df _{HK} (\xi )}{d\xi } \right] \left( \prod _{i=1}^3 (\wp(\xi )-e_i) ^{\tilde{l}_i/2} \right) . \nonumber 
\end{align}
If $ l_0 \neq 1/2$ (resp. $l_0 = 1/2$) and the integrand in Eq.(\ref{eq:intrepHeune}) (resp. Eq.(\ref{eq:intrepHeune00})) has a non-zero finite limit as $\mu \rightarrow 0$  (resp. $s \rightarrow 0$) for some $m$,
then the functions
\begin{align}
& f^{(i)} (x)=  \left( \prod _{j=1}^3 (\wp (x)-e_j )^{-l_j/2 } \right) \int _{-x+2\omega _i}^{x} \lim _{\mu \rightarrow 0} F(\xi; \lambda ,\mu , m) \frac{(\wp(x) -\wp (\xi ))^{\tilde{l}}}{(\tilde{\wp }(\xi )-\lambda )} d\xi \label{eq:intrepHeune}, \\
& \lambda = \frac{E+ (l_3-l_1)(2l_0+l_1+l_3) e_1+(l_3-l_2)(2l_0+l_2+l_3) e_2}{(e_1-e_2)(l_1+l_2+l_3+l_0)(2l_0-1)} +\frac{e_1}{e_2-e_1} , \quad (l_0 \neq 1/2), \nonumber \\
& f^{(i)} (x)=  \left( \prod _{j=1}^3 (\wp (x)-e_j )^{-l_j/2 } \right) \int _{-x+2\omega _i}^{x} \lim _{s \rightarrow 0} F(\xi; c/s ,bs^2 , m) \frac{(\wp(x) -\wp (\xi ))^{\tilde{l}}}{(\tilde{\wp }(\xi )-c/s )} d\xi \label{eq:intrepHeune00}, \\
& bc^2=  \frac{E+ (l_3-l_1)(l_1+l_3+1 ) e_1+(l_3-l_2)(l_2+l_3+1) e_2}{(e_1-e_2)(2l_1+2l_2+2l_3+1)} ,\quad (l_0=1/2), \nonumber 
\end{align}
$(i=0,1,2,3)$ are solutions to the elliptical representation of Heun's equation (Eq.(\ref{InoEF0})).
\end{thm}
We consider the limits $\lambda +\kappa _1/\mu \rightarrow 0,1,t$.
The following equations are obtained by setting $\lambda =-\kappa _1/\mu $, $\lambda =1-\kappa _1/\mu $, $\lambda =t-\kappa _1/\mu $ in the Fuchsian equation $D_{y_1}(\theta _0, \theta _1, \theta _t, \theta _{\infty}; \lambda ,\mu )$:
\begin{align}
& \frac{d^2y}{dz^2} + \left( \frac{-\tilde{\theta }_0}{z}+\frac{1-\tilde{\theta }_1}{z-1}+\frac{1-\tilde{\theta }_t}{z-t} \right) \frac{dy}{dz} +\frac{\tilde{\kappa }_1 (\tilde{\kappa }_2+1)z +t\tilde{\theta }_0\mu }{z(z-1)(z-t)}y=0, \label{Heun0} \\
& \frac{d^2y}{dz^2} + \left( \frac{1-\tilde{\theta }_0}{z}+\frac{-\tilde{\theta }_1}{z-1}+\frac{1-\tilde{\theta }_t}{z-t} \right) \frac{dy}{dz} +\frac{\tilde{\kappa }_1 (\tilde{\kappa }_2+1)(z-1) +(1-t)\tilde{\theta }_1\mu }{z(z-1)(z-t)}y=0, \label{Heun1} \\
& \frac{d^2y}{dz^2} + \left( \frac{1-\tilde{\theta }_0}{z}+\frac{1-\tilde{\theta }_1}{z-1}+\frac{-\tilde{\theta }_t}{z-t} \right) \frac{dy}{dz} +\frac{\tilde{\kappa }_1 (\tilde{\kappa }_2+1)(z-t) +t(t-1)\tilde{\theta }_t\mu }{z(z-1)(z-t)} y=0, \label{Heunt} 
\end{align}
where the parameters are defined as for the case $\mu =0$.
For the case $\tilde{\theta }_i =0$ $(i=0,1,t)$, we set $\mu=c/s$, $\lambda = i-\kappa _1/\mu  +bs^2$, and consider the the limit $s \rightarrow 0$.
Then we have
\begin{align}
& \frac{d^2y}{dz^2} + \left(\frac{1-\tilde{\theta }_1}{z-1}+\frac{1-\tilde{\theta }_t}{z-t} \right) \frac{dy}{dz} +\frac{\tilde{\kappa }_1 (\tilde{\kappa }_2+1)z -tbc^2 }{z(z-1)(z-t)}y=0, \label{Heun0-0} \\
& \frac{d^2y}{dz^2} + \left( \frac{1-\tilde{\theta }_0}{z}+\frac{1-\tilde{\theta }_t}{z-t} \right) \frac{dy}{dz} +\frac{\tilde{\kappa }_1 (\tilde{\kappa }_2+1)(z-1) +(t-1)bc^2}{z(z-1)(z-t)}y=0, \label{Heun1-0} \\
& \frac{d^2y}{dz^2} + \left( \frac{1-\tilde{\theta }_0}{z}+\frac{1-\tilde{\theta }_1}{z-1} \right) \frac{dy}{dz} +\frac{\tilde{\kappa }_1 (\tilde{\kappa }_2+1)(z-t) +t(1-t)bc^2 }{z(z-1)(z-t)} y=0. \label{Heunt-0} 
\end{align}
Note that we have similar propositions to Theorem \ref{thm:interepHeun}.

We consider the integral representations of solutions to Heun's equation for the case $\gamma =\delta =\epsilon =1$ and $\alpha = 3/2, \beta =1/2$, i.e. the case $l_0=1/2$, $l_1=l_2=l_3=-1/2$.
Recall that the functions
\begin{align}
& f_i(x)=\int_{-x+2\omega _i}^x \frac{e^{(\kappa +\zeta (\alpha ))\xi }\sigma (x)\sigma (\xi -\alpha )}{\sqrt{\sigma (x-\xi )\sigma (x+\xi )}}d\xi ,  \quad z=\frac{\wp (x)-e_1}{e_2-e_1} , \label{eq:intrep00010}
\end{align}
for $i=0,1,2,3$ are solutions to the Fuchsian differential equation $D_{y_1} (0,0,0,1;( \wp (\alpha ) -e_1)/(e_2-e_1) ,-(e_2-e_1 )\kappa /\wp ' (\alpha ))$ (see Proposition \ref{prop:intrep0001}).
The condition $s\rightarrow 0$ in Theorem \ref{thm:interepHeun} implies the condition $\alpha \rightarrow 0$ while setting $\kappa =-\zeta(\alpha )+\tilde{\kappa }$.
Therefore, it follows from Eq.(\ref{eq:intrep00010}) that the functions
\begin{align}
&  f_i(x)=\int_{-x+2\omega _i}^x \frac{e^{\tilde{\kappa }\xi }\sigma (x)\sigma (\xi  )}{\sqrt{\sigma (x-\xi )\sigma (x+\xi )}}d\xi , \label{eq:intrepH0001}
\end{align}
for $i=0,1,2,3$ are solutions to Heun's equation
\begin{align}
& \frac{d^2y}{dz^2} + \left( \frac{1}{z}+\frac{1}{z-1}+\frac{1}{z-t}\right) \frac{dy}{dz} +\frac{3z +(3e_1 -\tilde{\kappa }^2)/(e_2-e_1)}{4z(z-1)(z-t)} y=0, \label{eq:Heuns70001} 
\end{align}
by setting $z=(\wp (x)-e_1)/(e_2-e_1)$, and the functions 
\begin{align}
 f^{(i)}(x)&= \left( \prod _{i=1}^3 (\wp (x) -e_i) \right) ^{1/4}\int ^{x}_{-x+2\omega _i} \frac{e^{\tilde{\kappa }\xi }\sigma (x)\sigma (\xi  )}{\sqrt{\sigma (x-\xi )\sigma (x+\xi )}}d\xi \\
 & =\left( \frac{\sigma (x-\omega_1 )\sigma (x-\omega _2) \sigma (x-\omega _3)}{\sigma (x)} \right)^{1/2} \int ^{x}_{-x+2\omega _i} \frac{e^{\tilde{\kappa }\xi }\sigma (\xi  )}{\sqrt{\sigma (x-\xi )\sigma (x+\xi )}}d\xi  , \nonumber
\end{align}
for $i=0,1,2,3$ are solutions to Heun's equation in elliptical form for the case $l_0=1/2$, $l_1=l_2=l_3=-1/2$, 
\begin{align}
& \left(-\frac{d^2y}{dx^2} + \frac{3}{4} \wp(x) -\frac{1}{4} \sum_{i=1}^3 \wp (x+\omega_i ) +\tilde{\kappa }^2 \right)f(x)=0, \label{InoEF00001} 
\end{align}
The monodromy matrix of solutions to Eq.(\ref{eq:Heuns70001}) can be expressed in the form of those in section \ref{sec:intrepFs0000} by substituting $\kappa =-\zeta (\alpha )+\tilde{\kappa }$ and $\alpha =0$.
In fact, if $e^{2\omega _1 \tilde{\kappa }} \neq 1$ then the functions $f_0 (x) $ and $f_1(x)$ are linearly independent, and the monodromy matrices are written as
\begin{align}
& (f^{\gamma _0}_0(x) , f^{\gamma _0}_1(x))=  (f_0(x) , f_1(x))
\left(
\begin{array}{cc}
1 &  2(1-e^{2\omega _1 \tilde{\kappa }}) \\
0 & 1 
\end{array}
\right) , \label{eq:monodmHeun} \\
& (f^{\gamma _1}_0(x) , f^{\gamma _1}_1(x))=  (f_0(x) , f_1(x))
\left(
\begin{array}{cc}
1 &  0 \\
-2(1-e^{-2\omega _1 \tilde{\kappa }}) & 1 
\end{array}
\right) ,\nonumber \\
& (f^{\gamma _j}_0(x) , f^{\gamma _j}_1(x)) \quad \quad \quad \quad \quad \quad \quad \quad \quad \quad \quad \quad \quad \quad \quad \quad \quad \quad \quad \quad (j=2,3) \nonumber \\
& \quad \quad = (f_0(x) , f_1(x))
\left(
\begin{array}{cc}
1+2\frac{(e^{2\omega _1 \tilde{\kappa }} - e^{2\omega _j \tilde{\kappa }})(e^{2\omega _j \tilde{\kappa }}-1 )}{(e^{2\omega _1 \tilde{\kappa }}-1)e^{2\omega _j \tilde{\kappa }}} &  2\frac{(e^{2\omega _1 \tilde{\kappa }} - e^{2\omega _j \tilde{\kappa }})^2 }{(e^{2\omega _1 \tilde{\kappa }}-1)e^{2\omega _j \tilde{\kappa }}} \\
-2\frac{(e^{2\omega _1 \tilde{\kappa }} - 1)^2 }{(e^{2\omega _1 \tilde{\kappa }}-1)e^{2\omega _j \tilde{\kappa }}} & 1-2\frac{(e^{2\omega _1 \tilde{\kappa }} - e^{2\omega _j \tilde{\kappa }})(e^{2\omega _j \tilde{\kappa }}-1 )}{(e^{2\omega _1 \tilde{\kappa }}-1)e^{2\omega _j \tilde{\kappa }}}
\end{array}
\right) ,\nonumber 
\end{align}
which are obtained by analytic continuation of Eqs.(\ref{eq:f0f1gam0}, \ref{eq:f0f1gam1}, \ref{eq:f0f1gamj}) on the limit $\alpha \rightarrow 0$.
The monodromy matrices of solutions to Eq.(\ref{eq:Heuns70001}) are written as products of the monodromy matrices in Eq.(\ref{eq:monodmHeun}) and the scalar that is determined by the branching of $(\sigma (x-\omega _1 )\sigma (x-\omega _2 )\sigma (x-\omega _3 )/\sigma (x))^{1/2 }$.
If $\tilde{\kappa }=0$, then the integrals in Eq.(\ref{eq:intrepH0001}) are written as
\begin{align}
\int _{\infty }^z \frac{dw}{\sqrt{(w-z)(w-e_1)(w-e_2)(w-e_3)}}, \quad \int _{e_i }^z \frac{dw}{\sqrt{(w-z)(w-e_1)(w-e_2)(w-e_3)}}, 
\end{align}
for $i=1,2,3$ by setting $w=\wp (\xi)$ and $z=\wp (x)$.
These integrals coincide with the formula for the density function on root asymptotics of spectral polynomials for the Lame operator discovered by Borcea and Shapiro \cite{BS} (see also \cite{TakW}). 

The limits $\lambda +\kappa _1/\mu \rightarrow 0,1,t$ correspond respectively to the limits $\alpha \rightarrow e_1, e_2, e_3$.
The functions 
\begin{align}
&  f_{i'}(x)=\int_{-x+2\omega _{i'}}^x \frac{e^{(\kappa +\eta _i) \xi }\sigma (x)\sigma (\xi -\omega_i )}{\sqrt{\sigma (x-\xi )\sigma (x+\xi )}}d\xi , \quad (i'=0,1,2,3 ) \label{eq:intrepH0001123}
\end{align}
are solutions to the following Heun's equations; 
\begin{align}
& \frac{d^2y}{dz^2} + \left( \frac{1}{z-1}+\frac{1}{z-t}\right) \frac{dy}{dz} +\frac{z -\kappa ^2/(e_2-e_1)}{4z(z-1)(z-t)}y=0, \quad (i=1), \label{Heun00}  \\
& \frac{d^2y}{dz^2} + \left( \frac{1}{z}+\frac{1}{z-t}\right) \frac{dy}{dz} +\frac{z -1-\kappa ^2/(e_2-e_1)}{4z(z-1)(z-t)}y=0, \quad (i=2), \label{Heun11} \\
& \frac{d^2y}{dz^2} + \left( \frac{1}{z}+\frac{1}{z-1}\right) \frac{dy}{dz} +\frac{z -t-\kappa ^2/(e_2-e_1)}{4z(z-1)(z-t)} y=0, \quad (i=3), \label{Heuntt} 
\end{align}
by setting $z=(\wp (x)-e_1)/(e_2-e_1)$, and we have similar results for Heun's equations in elliptical form for the case $l_0=-l_1=l_2=l_3=-1/2$, $l_0=l_1=-l_2=l_3=-1/2$, $l_0=l_1=l_2=-l_3=-1/2$ respectively.
The monodromy matrices are expressed in similar forms as Eq.(\ref{eq:monodmHeun}).

{\bf Acknowledgments}
The author would like to thank Galina Filipuk and Yoshishige Haraoka for fruitful discussions and valuable comments.
Thanks are also due to Philip Boalch.
He is supported by the Grant-in-Aid for Young Scientists (B) (No. 19740089) from the Japan Society for the Promotion of Science.


\begin{thebibliography}{9999}
\bibitem{BE}
Belokolos E. D. and Enolskii V. Z., Reduction of Abelian functions and algebraically integrable systems. II, {\it J. Math. Sci. (New York)} {\bf 108} (2002), 295--374.
\bibitem{Boa}
Boalch P., From Klein to Painleve via Fourier, Laplace and Jimbo {\it Proc. London Math. Soc. (3)} {\bf 90} (2005), 167--208. 
\bibitem{BS}
Borcea J. and Shapiro B., Root asymptotics of spectral polynomials for the Lame operator, Preprint, math.CA/0701883.
\bibitem{DR1}
Dettweiler M. and  Reiter S., An algorithm of Katz and its application to the inverse Galois problem. Algorithmic methods in Galois theory {\it J. Symbolic Comput.} {\bf 30} (2000), 761--798.
\bibitem{DR2}
Dettweiler M. and Reiter S., On the middle convolution, Preprint, math.AG/0305311.
\bibitem{Fil}
Filipuk G., On the middle convolution and birational symmetries of the sixth Painleve equation, {\it Kumamoto J. Math.} {\bf 19} (2006), 15--23. 
\bibitem{GW}
Gesztesy F. and Weikard R., Treibich-Verdier potentials and the stationary (m)KdV hierarchy, {\it Math. Z.} {\bf 219} (1995), 451--476. 
\bibitem{Man}
Manin Yu. I., Sixth Painleve equation, universal elliptic curve, and mirror of $\bold P\sp 2$, Geometry of differential equations, 131--151, Amer. Math. Soc. Transl. Ser. 2, 186, Amer. Math. Soc., Providence, RI, 1998.
\bibitem{Oka}
Okamoto K., Studies on the Painleve equations. I. Sixth Painleve equation $P\sb {{\rm VI}}$, {\it Ann. Mat. Pura Appl. (4)}  {\bf 146} (1987), 337--381.
\bibitem{Ron}
Ronveaux A.(ed.), {\it Heun's differential equations,} Oxford Science Publications, Oxford University Press, Oxford, 1995.
\bibitem{Smi}
Smirnov A. O.,  Elliptic solitons and Heun's equation, {\it The Kowalevski property,} 287--305, CRM Proc. Lecture Notes, 32, Amer. Math. Soc., Providence (2002).
\bibitem{Tks}
Takasaki K., Painleve-Calogero correspondence revisited, {\it J. Math. Phys.} {\bf 42} (2001), 1443--1473.
\bibitem{Tak1}
Takemura K., The Heun equation and the Calogero-Moser-Sutherland system I: the Bethe Ansatz method. {\it Comm. Math. Phys.} {\bf 235} (2003), 467--494.
\bibitem{Tak2}
Takemura K., The Heun equation and the Calogero-Moser-Sutherland system II: the perturbation and the algebraic solution, {\it  Electron. J. Differential Equations} {\bf 2004} no. 15 (2004), 1--30. 
\bibitem{Tak3}
Takemura K., The Heun equation and the Calogero-Moser-Sutherland system III: the finite gap property and the monodromy, {\it J. Nonlinear Math. Phys.} {\bf 11} (2004), 21--46.
\bibitem{Tak4}
Takemura K., The Heun equation and the Calogero-Moser-Sutherland system IV: the Hermite-Krichever Ansatz, {\it Comm. Math. Phys.} {\bf 258} (2005), 367--403.
\bibitem{Tak5}
Takemura K.,  The Heun equation and the Calogero-Moser-Sutherland system V: generalized Darboux transformations, {\it J. Nonlinear Math. Phys.} {\bf 13} (2006), 584--611.
\bibitem{TakP}
Takemura K., Fuchsian equation, Hermite-Krichever Ansatz and Painlev\'e equation, Preprint, math.CA/0501428.
\bibitem{TakHP}
Takemura K., Heun equation and Painlev\'e equation, Preprint, math.CA/0503288.
\bibitem{TakW}
Takemura K., Finite-gap potential, Heun's differential equation and WKB analysis, Preprint, math.CA/0703256.
\bibitem{TV}
Treibich A. and Verdier J.-L., Revetements exceptionnels et sommes de 4 nombres triangulaires, {\it Duke Math. J.} {\bf 68} (1992), 217--236. 

\end{thebibliography}
\end{document}